\documentclass[11pt, leqno]{article}

\usepackage{amsmath,amssymb,amsthm, epsfig}

\usepackage{hyperref}

\usepackage{amsmath}

\usepackage[pdftex]{color}

\usepackage{mathtools}

\usepackage{color}

\usepackage{ulem}

\usepackage{mathrsfs}

\usepackage{wasysym}

\usepackage{stackrel}

\usepackage{pgfplots}
\pgfplotsset{compat=newest}

\title{\large{\bf Sharp and improved regularity estimates for weighted quasilinear elliptic equations of $p-$Laplacian type and applications}}
\author{\it by \smallskip \\ Jo\~{a}o Vitor  da Silva
\footnote{\noindent Universidade Estadual de Campinas - UNICAMP. Departamento  de Matemática. Campinas - SP, Brazil. \noindent \texttt{E-mail address: \url{jdasilva@unicamp.br}}}, \quad Disson S. dos Prazeres\footnote{\noindent Universidade Federal de Sergipe - UFS. Departamento  de Matemática, Brazil. \noindent \texttt{E-mail address: \url{disson@mat.ufs.br}}}, \quad Gleydson C. Ricarte\footnote{\noindent Universidade Federal do Cear\'{a} - UFC. Departamento de Matemática. Fortaleza - CE, Brazil. \noindent \texttt{E-mail address: \url{ricarte@mat.ufc.br}}} \\ $\&$\\ Ginaldo S. S\'{a}\footnote{\noindent Universidade Estadual de Campinas - UNICAMP. Departamento  de Matemática. Campinas - SP, Brazil. \noindent \texttt{E-mail address: \url{ginaldo@ime.unicamp.br}}}
}

\newlength{\hchng}
\newlength{\vchng}
\def \R {\mathbb{R}}

\def \div {\mathrm{div}}

\newcommand{\defeq}{\mathrel{\mathop:}=}

\newtheorem{theorem}{Theorem}[section]

\newtheorem{lemma}[theorem]{Lemma}

\newtheorem{proposition}[theorem]{Proposition}

\newtheorem{corollary}[theorem]{Corollary}

\theoremstyle{definition}

\newtheorem{definition}[theorem]{Definition}

\newtheorem{example}[theorem]{Example}

\newtheorem{remark}[theorem]{Remark}

\numberwithin{equation}{section}

\newcommand{\intav}[1]{\mathchoice {\mathop{\vrule width 6pt height 3 pt depth  -2.5pt
\kern -8pt \intop}\nolimits_{\kern -6pt#1}} {\mathop{\vrule width
5pt height 3  pt depth -2.6pt \kern -6pt \intop}\nolimits_{#1}}
{\mathop{\vrule width 5pt height 3 pt depth -2.6pt \kern -6pt
\intop}\nolimits_{#1}} {\mathop{\vrule width 5pt height 3 pt depth
-2.6pt \kern -6pt \intop}\nolimits_{#1}}}

\begin{document}
\maketitle

\begin{abstract}

In this manuscript, we obtain sharp and improved regularity estimates for weak solutions of weighted quasilinear elliptic models of Hardy-H\'{e}non-type, featuring an explicit regularity exponent depending only on universal parameters. Our approach is based on geometric tangential methods and uses a refined oscillation mechanism, compactness, and scaling techniques. In some specific scenarios, we establish higher regularity estimates and non-degeneracy properties, providing further geometric insights into such solutions. Our regularity estimates both enhance and, to some extent, extend the results arising from the $C^{p^{\prime}}$ conjecture for the $p$-Laplacian with a bounded source term. As applications of our results, we address some Liouville-type results for our class of equations. Finally, our results are noteworthy, even in the simplest model case governed by the $p$-Laplacian with regular coefficients:
$$
    \div\left( |\nabla u|^{p-2}\mathfrak{A}(|x|) \nabla u\right)  = |x|^{\alpha}u_+^m(x) \quad \text{in} \quad B_1
$$
under suitable assumptions on the data, with possibly singular weight $\mathfrak{h}(|x|) = |x|^{\alpha}$, which includes the Matukuma and Batt–
Faltenbacher–Horst's equations as toy models.

	\medskip
	\noindent \textbf{Keywords}: Hardy-H\'{e}non-type equations, quasilinear operators, sharp regularity estimates.
	\vspace{0.2cm}
	
	\noindent \textbf{AMS Subject Classification: Primary 35B65, 35J62; Secondary 35R35}
\end{abstract}

\newpage

\section{Introduction}

In this work, we establish sharp regularity estimates for weak solutions of quasilinear elliptic equations of $p$-Laplacian type with varying coefficients:
\begin{equation}\label{pobst}
	 \div \,\mathfrak{a}(x,\nabla u)  =  f(|x|, u(x)) \quad \text{in} \quad B_1,
\end{equation}
where $B_1 \subset \mathbb{R}^n$ denotes the unit $n$-dimensional ball with $n \geq 2$. Additionally, the vector field $\mathfrak{a}: B_1 \times \mathbb{R}^n \rightarrow \mathbb{R}^n$ is $C^1$-regular in the gradient variable and satisfies the following structural conditions (growth and ellipticity assumptions): for all $x, y \in B_1$ and $\xi, \eta \in \mathbb{R}^n$, we have:
\begin{equation} \label{condestr}
\left\{
\begin{array}{rclcl}
|\mathfrak{a}(x,\xi)|+|\partial_{\xi}\mathfrak{a}(x,\xi)||\xi| & \leq & \Lambda |\xi|^{p-1}, & & \\
\lambda |\xi|^{p-2}|\eta|^2 & \le & \langle \partial_{\xi}\mathfrak{a}(x,\xi)\eta,\eta \rangle, &  & \\
\displaystyle	\sup_{x, y \in B_1 \atop{x \ne y, \,\,\,|\xi|\ne 0 }} \frac{|\mathfrak{a}(x,\xi)-\mathfrak{a}(y,\xi)|}{\omega(|x-y|)|\xi|^{p-1}} & \le & \mathfrak{L}_0 < \infty, & &  \\
\end{array}
\right.
\end{equation}
where $1 < p < \infty$, $0 < \lambda \leq \Lambda < \infty$, $\mathfrak{L}_0 > 0$, and $\omega: [0,\infty) \rightarrow [0,\infty)$ represents a given modulus of continuity for the coefficients of the operator $\mathfrak{a}$, i.e., a concave and non-decreasing function with $\omega(0) = 0$.

It is well-known that the following Dini-type condition on the coefficients
\begin{equation}\label{EqDini-Cond}
\left\{
\begin{array}{rcl}
  \displaystyle \int_{0}^{R} \frac{\omega(t)}{t}\,dt < \infty & \text{for} & 1 < p < 2, \\
  \displaystyle \int_{0}^{R} \frac{\omega(t)^{\frac{2}{p}}}{t} \,dt < \infty &  \text{for} & p \geq 2,
\end{array}
\right.
\end{equation}
is essentially optimal for achieving $C^1$ estimates for quasilinear models with varying coefficients, as in the following model:
\begin{equation}\label{EqModel}
   -\div \,\mathfrak{a}(x,\nabla u) =  f(x) \quad  \text{in}  \quad \Omega, \quad \text{for} \quad f \in L(n, 1)(\Omega) \subset L^n(\Omega),
\end{equation}
where for $s \in (0, \infty)$ and $t \in (0, \infty]$, $L(s, t)(\Omega)$ denotes the Lorentz space. 

We recommend that the reader consult the seminal works of Mingione \textit{et al.}, such as \cite{DuzMing10} and \cite{KuusiMing12}, as well as the work of Dong and Zhu \cite{DZ23}, for enlightening contributions to this topic.

In this context, in this manuscript, we shall consider:
\begin{equation} \label{intdoscoef}
\omega(s) \lesssim s^{1+\sigma} \quad \text{and} \quad f \in \Xi(B_1 \times \mathfrak{I}) \,\, (\text{a suitable Banach space}),
\end{equation}
for $|\alpha| < \sigma$, where $\mathrm{c}_n > 0$ is a universal constant, and $f_0 \in \Xi(B_1)$. Moreover, we assume that
\begin{equation}\label{EqHomog-f}
|f(|rx|, st)| \leq \mathrm{c}_n r^{\alpha} s^m |f_0(x, t)| \quad \forall \,\, (x, t) \in B_1 \times \mathfrak{I}, \quad \text{and} \quad r, s \in (0, 1),
\end{equation}
for some $\alpha \in (-1-m, |p-2(1+m)|)$ and $0 \leq m < p - 1$. 

Additionally, we assume that
$$
\Xi(B_1\times \mathfrak{I}) =
\left\{
\begin{array}{lcl}
C^0(B_1 \times \mathfrak{I}) & \text{if} & 0 < \alpha < |p-2-m|, \\
L^{q}(B_1 \times \mathfrak{I}), \,\,  q > n & \text{if} & \max\left\{-1-m, -\frac{n}{q}\right\} < \alpha < 0.
\end{array}
\right.
$$

It is important to note that these conditions could be further relaxed. However, we choose to state our estimates based on condition \eqref{intdoscoef} to keep the presentation and estimates clearer.

\begin{example}
A typical model for \eqref{pobst} is the quasilinear elliptic Hardy-H\'{e}non-type problem governed by the $p$-Laplacian with varying coefficients:
$$
\displaystyle \div \left( |\nabla u|^{p-2} \mathfrak{A}(x) \nabla u \right) = \sum_{i=1}^{k_0} c_i |x|^{\alpha_i} u_+^{m_i}(x) \quad \text{in} \quad B_1,
$$
where $0 < \lambda \text{Id}_n \leq \mathfrak{A}(\cdot) \leq \Lambda \text{Id}_n$ is a symmetric matrix with sufficiently smooth entries (cf. \cite{pucciserrin}), $k_0 \in \mathbb{N}$, $c_i \in \mathbb{R}_+$, and
$$
0 < \alpha_i < |p-2(1+m)| \quad \text{and} \quad 0 \leq m_i < p - 1 \quad (\text{or} \quad -1 < \alpha_i < 0) \quad \text{for} \quad 1 \leq i \leq k_0.
$$
In this case, $f(x, u)$ satisfies \eqref{EqHomog-f} as follows:
$$
\displaystyle f(rx, su) \leq r^{\displaystyle \min_{1 \leq i \leq k_0} \{\alpha_i\}} s^{\displaystyle \min_{1 \leq i \leq k_0} \{m_i\}} f(x, u) \quad \forall \,\, (x, u) \in B_1 \times \mathfrak{I},
$$
$$
\qquad \left( \text{resp.} \,\,\, f(rx, su) \leq r^{\displaystyle \max_{1 \leq i \leq k_0} \{\alpha_i\}} s^{\displaystyle \min_{1 \leq i \leq k_0} \{m_i\}} f(x, u) \right).
$$
where
$$
\displaystyle \|f_0\|_{L^\infty(B_1)}  := \sup_{ \overline{B_1}} f(|x|, u(x))
$$
\end{example}

Before stating our main results, we introduce some notation and classes of functions useful for our purposes.

For subsequent analysis, we define the class:
$$
\mathbb{H}_{\Lambda, \lambda, n} \defeq
\left\{
\begin{array}{ccc}
\mathfrak{h} \in W^{1, p}(B_1) & | & \div \, \mathfrak{a}(\nabla \mathfrak{h}) = 0 \quad \text{in} \quad B_1 \,\, \text{in the weak sense}, \\
& & \text{for a vector field} \,\, \mathfrak{a}: \mathbb{R}^n \to \mathbb{R}^n \,\, \text{satisfying} \,\, \eqref{condestr}
\end{array}
\right\}.
$$

We note that, although this constitutes a large class of functions, it is well established (see Section \ref{Motivat} for a detailed presentation) that there exists a universal modulus of continuity for the gradient of profiles in \( \mathbb{H}_{\Lambda, \lambda, n} \). Specifically, if \( \mathfrak{h} \in \mathbb{H}_{\Lambda, \lambda, n} \), then there exist constants \( \mathrm{C}^{\star} > 0 \) and \( \alpha^{\star} \in (0, 1) \), depending only on the dimension and structural parameters, such that
\begin{equation}\label{EqHomogeneous}
  \|\mathfrak{h}\|_{C^{1,\alpha^{\star}}(B_{1/2})} \leq \mathrm{C}^{\star} \cdot \|\mathfrak{h}\|_{L^{\infty}(B_1)}.
\end{equation}

Additionally, it is important to highlight that the exponent \( \alpha^{\star} \in (0, 1) \) represents a theoretical bound in the regularity theory of weak solutions of
$$
\div \, \mathfrak{a}(\nabla \mathfrak{h}) = 0 \quad \text{in} \quad B_1,
$$
which can be characterized as follows:
\begin{equation}\label{Def_alpha_Hom}
\alpha_{\mathrm{H}} \defeq
\sup \left\{
\begin{array}{ccc}
\alpha^{\star} \in (0, 1) & | & \text{there exists} \, \mathrm{C}_{n, \alpha^{\star}} > 0 \, \text{such that} \, \eqref{EqHomogeneous} \\
& &  \text{holds with} \, \mathrm{C}_{n, \alpha^{\star}} = \mathrm{C}^{\star} \, \text{for all} \, \mathfrak{h} \in \mathbb{H}_{\Lambda, \lambda, n}
\end{array}
\right\}.
\end{equation}

\subsection*{Main Theorems}

\begin{theorem}[{\bf Sharp/improved regularity estimates}]\label{sharp_regularity}
Let \( u \in W^{1,p}(B_1) \) be a bounded weak solution to \eqref{pobst} for \( f(|x|, u) = f(|x|) \). Suppose further that the assumptions \eqref{condestr} and \eqref{intdoscoef} hold. Then, \( u \in C^{1, \beta}_{loc}(B_1) \), where
\begin{equation}\label{sharp_exponent}
\beta = \min \left\{ \alpha_{\mathrm{H}}^-, (\alpha + 1) \cdot \min \left\{ 1, \frac{1}{p-1} \right\} \right\},
\end{equation}
and \( \alpha_{\mathrm{H}} \in (0, 1) \) comes from \eqref{Def_alpha_Hom}.

Moreover, the following estimates hold:
$$
\|u\|_{C^{1, \beta}(B_{1/2})} \leq
\left\{
\begin{array}{ccl}
  \mathrm{C}_{p, n, \alpha} \left(\|u\|_{L^{\infty}(B_1)} + \|f\|^{\frac{1}{p-1}}_{L^{\infty}(B_1)}\right) & \text{if} & 0 < \alpha < |p-2|, \\
  \mathrm{C}_{p, n, \alpha} \left(\|u\|_{L^{\infty}(B_1)} + \|f\|^{\frac{1}{p-1}}_{L^{q}(B_1)}\right) & \text{if} & -1 < \alpha < 0 \quad \text{and} \quad p > 2.
\end{array}
\right.
$$
\end{theorem}

It is important to note that our approach is a byproduct of a new oscillation-type
estimate, initially considered in \cite{LL} and \cite{ATU17}, combined with a
localized analysis. The proof is carried out by studying two complementary cases:
\begin{itemize}
    \item If the gradient is small with a controlled magnitude (i.e., \( |\nabla u| \lesssim \text{O}(r^{\beta}) \)), then a balanced perturbation of the \( p \)-harmonic profile leads to the inhomogeneous problem at the limit, via a stability argument in the \( C^1 \)-fashion.
    \item Conversely, if the gradient has a uniform lower bound, i.e., \( |Du| \geq \mathfrak{L}_{\ast} > 0 \), then classical estimates (see, e.g., \cite[Theorem 3.13]{QhanFlin}) can be applied since the operator becomes uniformly elliptic:
    \[
        \text{div}(\mathfrak{A}^{\ast}(x)\nabla u) \lesssim  \mathrm{C}_0(\mathfrak{L}_{\ast}^{-1}, \| f \|_{\Xi(B_1)}).
    \]
\end{itemize}

\begin{example}\label{Examp01}
Let \( 1< p < \infty \) and \( \alpha \in (0, |p-2|) \) be fixed constants. Then, consider the function
\[
u(x_1, \ldots, x_n) =  \frac{|x_i|^{1+\frac{1+\alpha}{p-1}}}{(1+\alpha)^{\frac{1}{p-1}}\left(1+\frac{1+\alpha}{p-1}\right)} \quad \text{for} \quad x \in B_1,
\]
which satisfies the following equation in the weak sense:
\[
\Delta_p u(x) = |x_i|^{\alpha} \lesssim |x|^{\alpha} \quad \text{in} \quad B_1.
\]
The singular (or critical) set is given by
\[
\mathcal{S}_u(B_1) = \{x \in B_1: |\nabla u| = 0\} = \{x \in B_1: x_i = 0\} \quad \text{with} \quad \mathscr{H}_{\text{dim}}(\mathcal{S}_u(B_1))= 1.
\]
In particular, we have that \( u \) belongs to \( C_{\text{loc}}^{1, \frac{1+\alpha}{p-1}}(B_1) \).
\end{example}


As a consequence of Theorem \ref{sharp_regularity}, we obtain sharp regularity estimates in specific scenarios.

\begin{corollary}\label{Corol01}
Let \( B_1 \subset \mathbb{R}^n \) and let \( u \in W^{1,p}(B_1) \) be a weak solution of \eqref{pobst} for \( p > 2 \) where \( f(|x|) \simeq |x|^{\alpha} \) and \( \alpha \in \left(0, p-2\right) \). Suppose further that \( \alpha_{\mathrm{H}} \in \left(\frac{1+\alpha}{p-1}, 1\right] \). Then, \( u \in C^{1,  \frac{\alpha+1}{p-1}}(B_{1/2}) \) and
\[
\|u\|_{C^{1, \frac{\alpha+1}{p-1}}(B_{1/2})} \leq \mathrm{C}_{p, n, \alpha} \left( \|u\|_{L^{\infty}(B_1)} + \|f\|^{\frac{1}{p-1}}_{L^\infty(B_1)}\right).
\]
\end{corollary}

Additionally, we can present an improved estimate for weak solutions in the plane (cf. \cite{ATU17}).

\begin{corollary}[{\bf Optimal estimates in $2$-$\mathrm{D}$}]\label{Corol02}
Let \( B_1 \subset \mathbb{R}^2 \) and let \( u \in W^{1,p}(B_1) \) be a bounded weak solution of
\[
-\Delta_p u = f(|x|), \quad p > 2,
\]
with \( f(|x|) \simeq |x|^{\alpha} \) and \( \alpha \in \left(0, \tau_0(p-1)\right) \) for a constant \( \tau_0 < \frac{p-2}{p-1} \) (see Proposition \ref{PropB-K}). Then, \( u \in C^{1, \frac{1+\alpha}{p-1}}(B_{1/2}) \) and
\[
\|u\|_{C^{1, \frac{1+\alpha}{p-1}}(B_{1/2})} \leq \mathrm{C}_{p, \alpha} \left( \|u\|_{L^{\infty}(B_1)} + \|f\|^{\frac{1}{p-1}}_{L^\infty(B_1)}\right).
\]
\end{corollary}

The previous regularity results (Corollaries \ref{Corol01} and \ref{Corol02}) can be interpreted as a type of ``Schauder estimate'' for weak solutions of the \( p \)-Laplacian. In this context, we achieve an improvement in regularity gain concerning the appropriate oscillation and homogeneity of the forcing term (cf. \cite{BCDS22} for related results on the \( \mathbf{p} \)-Laplace system).

To the best of our knowledge, this aspect has not yet been addressed in the classical literature and plays a crucial role in the regularity theory of nonlinear models. Thus, it contributes to advancing research in the theory of quasilinear PDEs and related topics, such as free boundary problems.

\begin{remark}
Let \( x_0 \in \Omega \) (with \( 0 \in \Omega \)), and consider the linear transformation
$$
\mathrm{T}_{x_0}: \mathbb{R}^n \rightarrow \mathbb{R}^n \quad \text{defined by} \quad \mathrm{T}_{x_0}(x) = x - x_0.
$$
If the model equation
\begin{equation}\label{EqWeight}
\mathrm{div} \mathfrak{a}(x, \nabla u(x)) = \mathrm{f}(|x|) \quad \text{in} \quad B_1
\end{equation}
is invariant under the transformation \( \mathrm{T}_{x_0} \), with
$$
\mathrm{f} \in C^0([0, \infty)) \quad \text{and} \quad |\mathrm{f}(|x|)| \leq \mathrm{c}_n |x|^{\alpha}, \quad \text{for} \quad \alpha > 0,
$$
then, according to Theorem \ref{sharp_regularity}, a weak solution belongs to \( C^{1, \beta} \) at \( x_0 \).

In particular, this holds when the weight is a multiple of the distance to a closed set: let \( \mathrm{F} \subset \subset \Omega \) be a fixed closed set, and consider \( \mathrm{f}(|x|) = \mathrm{dist}^{\alpha}(x, \mathrm{F}) \) for a given \( \alpha > 0 \). In this case, weak solutions to \eqref{EqWeight} belong to \( C^{1, \beta} \) along \( \mathrm{F} \). Furthermore,
$$
\sup_{B_r(x_0) \atop{x_0 \in \mathrm{F}}} |u(x) - u(x_0) - \nabla u(x) \cdot (x - x_0)| \leq \mathrm{C} \cdot r^{1+\beta},
$$
where \( \mathrm{C}>0 \) is a universal constant.
\end{remark}

In our second main result, we obtain sharp higher growth estimates for weak solutions to \( p \)-Laplacian-type operators \eqref{pobst} at critical points of existing solutions when \( \alpha > -1 - m \). 

\begin{theorem}[{\bf Higher Regularity Estimates}]\label{Hessian_continuity}
Suppose the assumptions of Theorem \ref{sharp_regularity} are satisfied. Additionally, assume that \( x_0 \in B_1 \) is a local extremum point for \( u \) and that \( f(|x|, t) \simeq |x - x_0|^{\alpha} t_+^m \) with \( \alpha + 1 + m > 0 \). Then,
\begin{equation}\label{Higher Reg}
\sup_{x \in B_r(x_0)} u(x) \leq \mathrm{C} r^{1 + \frac{1 + \alpha + m}{p - 1 - m}},
\end{equation}
for \( r \in (0, 1/2) \) and \( \mathrm{C} > 0 \) being a universal constant.
\end{theorem}

\medskip

\begin{remark}
We would like to highlight some consequences of the previous Theorem \ref{Hessian_continuity}. Let \( \beta = \frac{1 + \alpha + m}{p - 1 - m} \). We observe that at critical points, the following relationships hold:

$$
\left\{
\begin{array}{ccccc}
  \beta \in (0, 1) & \Leftrightarrow & \alpha < p - 2(1 + m) & \Rightarrow & u \in C^{1, \beta} \text{ at } x_0,\\
  \beta = 1 & \Leftrightarrow & \alpha = p - 2(1 + m) & \Rightarrow & u \in C^{1, 1} \text{ at } x_0, \\
  \beta > 1 & \Leftrightarrow & \alpha > p - 2(1 + m) & \Rightarrow & u \in C^{2} \text{ at } x_0.
\end{array}
\right.
$$

In particular, weak solutions enjoy classical estimates along the set of critical points provided \( \alpha > p - 2(1 + m) \) (see, e.g., Example \ref{Examp01}).
\end{remark}

\begin{remark}[{\bf Equations of Matukuma and Batt–
Faltenbacher–Horst type}]\label{Remark-M-BFH}

We also study the related elliptic equation
\begin{equation}\label{EqMatukuma}
\mathrm{div}\left(\mathfrak{a}(|x|) |\nabla u|^{p-2} \nabla u\right) = \mathfrak{h}(|x|) f(u) \quad \text{in } \quad \Omega \subset \mathbb{R}^n, \quad  p > 1,
\end{equation}
where \(\mathfrak{a}, \mathfrak{h} : \mathbb{R}^+ \to \mathbb{R}^+\) are radial functions of class \(C^1(\mathbb{R}^+)\) and $C^0(\mathbb{R}^+)$, respectively. By way of motivation, the celebrated Matukuma equation (resp. Batt–
Faltenbacher–Horst's equation, see section \ref{Subsection1.2} for more details) is a prototype for \eqref{EqMatukuma}. Moreover, the strong absorption term satisfies

\begin{itemize}
    \item[\textbf{(F1)}] $f \in C^0(\mathbb{R})$;
    \item[\textbf{(F2)}] $f(0) = 0$, $f$ is nondecreasing on $\mathbb{R}$, and $f(t) > 0$ for $t > 0$.
\end{itemize}

Now, consider the model equation
\begin{equation}\label{ModelEq}
    \mathrm{div} \left( |x|^k | \nabla u|^{p-2} \nabla u \right) = |x|^{\alpha} f(u)
\end{equation}
Hence, if $f(u) = \mathrm{C}u_{+}^{m}$, for some $\mathrm{C} > 0$ to be determined later, straightforward computation shows that the function
\[
u(x) = |x|^{\beta}, \quad \text{for} \quad \beta = \frac{\alpha + p - k}{p - 1 - m},
\]
is a nonnegative entire radial solution of \eqref{ModelEq}, provided that
\[
\mathrm{C} = \mathrm{C}_{n, p, \beta, k} = \beta^{p-1} \left[\beta(p - 1) + k + n - p\right].
\]
Observe that $\mathrm{C} > 0$, and $u$ is also of class $C^{\infty}(\mathbb{R}^n \setminus \{0\})$, but only H\"{o}lder continuous at the origin if $\beta \in (0, 1)$, that is, whenever $m + 1 < k - \alpha < p$.

In this context, and motivated by Theorem \ref{Hessian_continuity}, we can address the following estimate for weak solutions to \eqref{ModelEq}:
    $$
\sup_{x \in B_r(x_0)} u(x) \leq \mathrm{C} r^{1 + \frac{1 + \alpha + m - k}{p - 1 - m}},
    $$
    for a universal constant $\mathrm{C}>0$, provided that \( x_0 \in B_1 \) is a free boundary point for \( u \), and \( f(|x|, t) \simeq |x - x_0|^{\alpha} t_+^m \), with \( \alpha + 1 + m > k \).
\end{remark}

\medskip

As a consequence of Theorem \ref{Hessian_continuity}, we also derive the sharp gradient growth rate at local extremum points.

\begin{corollary}[{\bf Sharp Gradient Growth}]\label{Ric}
Under the assumptions of Theorem \ref{Hessian_continuity}, there exists a universal constant \( \mathrm{C} > 0 \) such that
\[
\sup_{B_r(x_0)} |\nabla u(x)| \leq \mathrm{C} \cdot \left(\|u\|_{L^\infty(\Omega)} + \|f\|_{\Xi(\Omega)}^{\frac{1}{p+1}}\right) r^\beta,
\]
for any point \( x_0 \in \partial\{u>0\} \) and for all \( 0 < r < \frac{1}{4} \).
\end{corollary}

\medskip

From now on, we will denote the critical zone of existing solutions as
\[
\mathcal{S}_{r,\beta} (u,\Omega) := \{x_0 \in \Omega \, : \, |\nabla u(x_0)| \leq r^\beta\}, \quad \text{for } 0 \leq r \leq 1.
\]
A geometric interpretation of Theorem \ref{sharp_regularity} states that if \( u \) solves \eqref{pobst} (with \( f(|x|, t) = f(|x|) \)) and \( x_0 \in \mathcal{S}_{r,\beta} (u,\Omega) \), then in a neighborhood of \( x_0 \), we have
\[
\sup_{B_r(x_0)} |u(x)| \leq |u(x_0)| + \mathrm{C} \cdot r^{1+\beta}.
\]
On the other hand, from a geometric viewpoint, it is crucial to obtain the corresponding sharp lower bound estimate for such operators with inhomogeneous degeneracy, providing a sort of qualitative information. This feature is termed the \textit{Non-degeneracy Property of Solutions}.

Next, we will consider the additional hypotheses
\begin{equation}\label{7.1}
   \sum_{i=1}^n\left| \dfrac{\partial \mathfrak{a}_i}{\partial x_i}(x,\xi)\right| \leq \Lambda_0 |\xi|^{p-1}  \qquad (\text{for every } \,\,\, \xi \in \mathbb{R}^n) 
\end{equation}
for some \( 0 \leq \Lambda_0 < \infty \).


Last but not least, in our final result, we will assume the condition \eqref{7.1} holds for all \( x \in B_1 \). Furthermore, we will assume that there exists a positive constant \( \mathrm{c}_0 \) such that the source term \( f \) satisfies a lower growth condition given by
\begin{equation}\label{7.2}
	\mathrm{c}_0 |x|^{\alpha} \leqslant f(|x|), \quad \text{a.e. in }\,\, B_1.
\end{equation}

\medskip
\begin{theorem}[{\bf Non-degeneracy Along Extremum Points}]\label{NãoDeg}
Assume the conditions \eqref{7.1} and \eqref{7.2} are satisfied. Let \( p \in (1, \infty) \) and \( u \in W^{1, p}(B_1) \) be a weak solution of the problem \eqref{pobst}. Then, there exists \( r^{\ast} > 0 \) such that for every local extremum point \( x_0 \in B_1 \) and every \( r \in (0, r^{\ast}) \) fulfilling \( B_r(x_0) \subset B_1 \), we have
\[
\sup_{\partial B_r(x_0)} \left( u(x) - u(x_0) \right) \geq \mathrm{C}_{\sharp}(n, p, \alpha, \Lambda, \Lambda_0)\cdot r^{1 + \frac{1+\alpha}{p-1}}
\]
for a universal constant $\mathrm{C}_{\sharp}>0$.
\end{theorem}

\begin{remark}\label{Rem_No-Deg}
We emphasize that a similar non-degeneracy result holds for a general forcing term satisfying the lower growth condition
\[
f(|x|, u) \geq \mathrm{c}_0 |x|^{\alpha} u_+^m(x) \quad \text{a.e. in} \quad B_1.
\]
In this context, it is possible to obtain a non-degeneracy estimate at critical points \( x_0 \in B_1 \) as follows:
\[
\sup_{x \in \partial B_r(x_0)} u(x) \geq \mathrm{C}_{\sharp}(n, p, m, \alpha, \mathrm{c}_0, \Lambda, \Lambda_0) \, r^{1+\frac{1+\alpha+m}{p-1-m}},
\]
for some explicit universal constant \( \mathrm{C}_{\sharp} > 0 \).
\end{remark}

\begin{example}[\textbf{Radial Scenario}]\label{RadialEx}
Under the assumptions from Example \ref{Examp01}, we will consider the radial function \( u(r) = \mathrm{c} r^\beta \), where \(\beta = 1 + \frac{1 + \alpha + m}{p - 1 - m}\) and \(r = |x|\).

Note that
\[
u^{\prime}(r) = \mathrm{c} \beta r^{\beta - 1}.
\]
Recall that the \( p \)-Laplacian for a radial profile is given by:
\[
\Delta_p u = \frac{d}{dr}\left(|u^{\prime}(r)|^{p-2}u^{\prime}(r)\right) + \left( \frac{n-1}{r} \right)|u^{\prime}(r)|^{p-2}u^{\prime}(r).
\]
First, observe that
\begin{itemize}
\item The term \( |u^{\prime}(r)|^{p-2} \) is:
\[
|u^{\prime}(r)|^{p-2} = \left| \mathrm{c} \beta r^{\beta - 1} \right|^{p - 2} = \mathrm{c}^{p - 2} \beta^{p - 2} r^{(\beta - 1)(p - 2)}.
\]
Then,
\begin{align*}
\frac{d}{dr}\left(|u^{\prime}(r)|^{p-2}u^{\prime}(r)\right) &= \frac{d}{dr}\left[\left(\mathrm{c}^{p - 2}  \beta^{p - 2} r^{(\beta - 1)(p - 2)}\right) \left(\mathrm{c} \beta r^{\beta - 1}\right)\right] \\[0.2cm]
&= \mathrm{c}^{p - 1} \beta^{p - 1} \frac{d}{dr} r^{(\beta - 1)(p - 1)} \\[0.2cm]
&= \mathrm{c}^{p - 1} \beta^{p - 1} (\beta - 1)(p - 1) r^{(\beta - 1)(p - 1) - 1}.
\end{align*}
\item The term \( \frac{n - 1}{r} u^{\prime}(r) \) is:
\[
\frac{n - 1}{r} u^{\prime}(r) = \mathrm{c} (n - 1) \beta r^{\beta - 1}.
\]
\end{itemize}
Now, we substitute these expressions into the \( p \)-Laplacian formula:
\begin{align*}
\Delta_p u  & = \mathrm{c}^{p - 1} \beta^{p - 1} (\beta - 1)(p - 1) r^{(\beta - 1)(p - 1) - 1} + \mathrm{c}^{p - 1} \beta^{p - 1} (n - 1) r^{(\beta - 1)(p - 1) - 1} \\
&= \mathrm{c}^{p - 1} \beta^{p - 1} \left[(\beta - 1)(p - 1) + (n - 1)\right] r^{(\beta - 1)(p - 1) - 1} \\
&= \mathrm{c}^{p - 1} \left(\frac{p + \alpha}{p - 1 - m}\right)^{p - 1} \left[\frac{(1 + \alpha + m)(p - 1) + (n - 1)(p - 1 - m)}{p - 1 - m}\right] r^{(\beta - 1)(p - 1) - 1}.
\end{align*}
Hence, if we choose
\[
\mathrm{c} = \left(\frac{p - 1 - m}{(1 + \alpha + m)(p - 1) + (n - 1)(p - 1 - m)}\right)^{\frac{1}{p - 1 - m}} \left(\frac{p - 1 - m}{p + \alpha}\right)^{\frac{p - 1}{p - 1 - m}},
\]
we obtain in the weak sense
\[
\Delta_p u(x) = |x|^\alpha u_+^m(x),
\]
thereby illustrating the universal constants from Theorem \ref{NãoDeg} and Remark \ref{Rem_No-Deg}.
\end{example}

As a result of Theorems \ref{Hessian_continuity} and \ref{NãoDeg}, we obtain a positive density result for the non-coincidence set.

\begin{corollary}\label{positive density}
Let \( u \) be a weak solution of \eqref{pobst} and let \( x_0 \in \{u > 0\} \cap B_{1/2} \). There exists a universal constant \( \mathrm{c}_0 > 0 \) such that
$$
\frac{\mathcal{L}^n(\{u > 0\} \cap B_r(x_0))}{\mathcal{L}^n(B_r(x_0))} \geq \mathrm{c}_0.
$$
\end{corollary}

As another consequence of Theorems \ref{Hessian_continuity} and \ref{NãoDeg}, we establish the porosity of the zero-level set.

\begin{definition}
A set \( \mathcal{S} \subset \mathbb{R}^n \) is called \textit{porous} with porosity \( \delta_{\mathcal{S}} > 0 \) if there exists \( R > 0 \) such that
\[
\forall x \in \mathcal{S}, \; \forall r \in (0, R), \; \exists y \in \mathbb{R}^n \; \text{such that} \; B_{\delta_{\mathcal{S}} r}(y) \subset B_r(x) \setminus \mathcal{S}.
\]
\end{definition}

We note that a porous set with porosity \( \delta_{\mathcal{S}} > 0 \) satisfies
$$
\mathscr{H}_{\mathrm{dim}}(\mathcal{S}) \leq n - c \delta_{\mathcal{S}}^n,
$$
where \( \mathscr{H}_{\mathrm{dim}} \) denotes the Hausdorff dimension and \( \mathrm{c} = \mathrm{c}(n) > 0 \) is a dimensional constant. In particular, a porous set has Lebesgue measure zero (see, for example, \cite{Zaj87}).

\begin{corollary}\label{Hausdorff dimension}
Let \( u \) be a non-negative weak solution of \eqref{pobst}. Then, the level set \( \partial \{ u > 0 \} \cap B_r \) is porous with a porosity constant that is independent of \( r > 0 \).
\end{corollary}

The proofs of Corollaries \ref{positive density} and \ref{Hausdorff dimension} are similar to those found in \cite{daSRosSal19} and \cite{daSSal18}. Therefore, we have chosen not to present the demonstrations here to avoid unnecessary repetition of the same arguments.

\subsection*{Some applications}

As an application of our findings, it is interesting to present a Liouville-type result for a class of quasilinear Hen\'{o}n-type equations. Specifically, we have:

\begin{theorem}[{\bf Liouville-type result $\mathrm{I}$}]\label{Liouville I}
Let \( u \) be a non-negative entire weak solution to
$$
\div \, \mathfrak{a} (x, \nabla u) = |x|^{\alpha}u^m(x) \quad \text{in} \quad \mathbb{R}^n,
$$
with \( u(0) = 0 \). Suppose that
\begin{equation}\label{Liouville_eq1}
\lim_{|x| \to \infty} \frac{u(x)}{|x|^{1+ \frac{1+\alpha+m}{p-1-m}}} = 0.
\end{equation}
Then, \( u \equiv 0 \).
\end{theorem}

The following Liouville-type result is a refined version of the previous theorem. For the sake of clarity, we will state it for the $p$-Laplacian operator.

\begin{theorem}[{\bf Liouville-type result $\mathrm{II}$}]\label{Liouville II}
Let \( u \) be a non-negative entire solution to
$$
\Delta_p u  = |x|^{\alpha}u^m(x) \quad \text{in} \quad \mathbb{R}^n.
$$
Suppose that
\begin{equation}\label{H-Liuville II}
\limsup_{|x| \to \infty} \frac{u(x)}{|x|^{1+ \frac{1+\alpha+m}{p-1-m}}} < \tau(n,p,m,\alpha),
\end{equation}
where
$$
\tau(n,p,m,\alpha) \coloneqq \left(\frac{p-1-m}{(1+\alpha+m)(p-1)+(n-1)(p-1-m)}\right)^\frac{1}{p-1-m} \left(\frac{p-1-m}{p+\alpha}\right)^{\frac{p-1}{p-1-m}},
$$
then \( u \equiv 0 \).
\end{theorem}

\begin{remark}
We would like to emphasize that Theorem \ref{Liouville II} represents a quantitative version of a Liouville-type result for quasilinear Hardy-H\'{e}non-type problems, while Theorem \ref{Liouville I} is a qualitative result. For this reason, Theorem \ref{Liouville II} implies the validity of Theorem \ref{Liouville I} (in the case of $p-$Laplacian).
\end{remark}

\subsection*{The borderline scenario}

In this section, we will turn our attention to the critical equation obtained as $m \to (p-1)^{+}$, i.e.,
\begin{equation}\label{Crit_Eq}
\Delta_p\,u = |x|^{\alpha}u^{p-1}(x) \quad \text{in} \quad \Omega.
\end{equation}

Observe that such an equation is considered critical, as all the regularity estimates established so far degenerate when $m \to p-1$. Moreover, we can rewrite equation \eqref{Crit_Eq} as
$$
\Delta_p\, u = |x|^{\alpha}u^{\kappa} u^{p-1-\kappa}(x)
$$
for any $\kappa>0$. In particular, we obtain from Theorem \ref{Hessian_continuity} that if $x_0 \in B_1$ satisfies $u(x_0) = 0$, then $D^{|\tau|} u(x_0) =0$, for all multi-indices $\tau \in \mathbb{N}^n$. Thus, any vanishing point is an infinite order zero. Finally, under the assumption that $u$ is an analytic function, one could conclude that $u \equiv 0$.

In this context, utilizing barrier arguments that leverage the scalar properties of the model equation \eqref{Crit_Eq}, we will demonstrate that a positive solution to \eqref{Crit_Eq} is indeed precluded from vanishing at interior points.

\begin{theorem}
Let $u \in W^{1, p}(\Omega)$ be a non-negative weak solution to \eqref{Crit_Eq}. Suppose that there exists an interior vanishing point $x_0 \in \Omega$, i.e., $u(x_0) =0$, then $u \equiv 0$ in $\Omega$.
\end{theorem}

\begin{proof}

To establish the desired result, we have opted to approximate equation \eqref{Crit_Eq} with a regularized problem. In effect, for $\epsilon > 0$ fixed , let $u_{\epsilon}$ be the weak solution to
$$
 \left\{
\begin{array}{rclcl}
\textrm{div}\,\mathfrak{a}(x,\nabla u_{\epsilon}) & = & \zeta_{\epsilon}(x)\cdot u_{\epsilon}^{p-1}(x) & \mbox{in} & \Omega \\
u_{\epsilon}(x) & = & \displaystyle\inf_{\partial \Omega} u(x) & \mbox{on} & \partial \Omega,
\end{array}
\right.
$$
where $\zeta_{\epsilon}(x) := (|x|+\epsilon)^{\alpha}$.

Now, our goal will be to analyze the limiting behavior as $\epsilon \to 0^{+}$. First, observe that the sequence of approximating functions $\{u_{\epsilon}\}_{\epsilon>0}$ converges locally to the weak solution $u$ of the original problem \eqref{Crit_Eq}. Moreover, note also that $\zeta_0(x) \le \zeta_{\epsilon}(x)$. Therefore, by the Comparison Principle (see Theorem \ref{comp}), we conclude that $u_{\epsilon} \leq u$ in $\Omega$.

Next, assume that there exists an interior point $x_0 \in \Omega$ such that $u(x_0)=0$. It follows that $u_{\epsilon}(x_0) = 0$ for all $\epsilon>0$. Since $\zeta_{\epsilon}(x)$ is non-degenerate, by \cite[Theorem 3.3]{daSSal18}, we have $u_{\epsilon} \equiv 0$ in $\Omega$. Thus, since $u_{\epsilon}$ converges to $u$ as $\epsilon \to 0$, it follows that $u \equiv 0$ in $\Omega$, thereby completing the proof.

\end{proof}

\subsection{State-of-the-Art and brief motivations of the model equation }\label{Motivat}

\subsection*{Hardy-H\'{e}non's equations: from astrophysics to diffusion phenomena}

The beginning of studies on H\'{e}non's equations dates back to the early 1970s, when Michel H\'{e}non, in \cite{Henon73}, used the concentric shell model to numerically explore the stability of spherical stellar systems in equilibrium under spherical perturbations, and examined a type of semilinear problem. In this investigation, it is inferred that the existence of stationary stellar dynamics models, see \cite{Batt77}, \cite{BFH86} and \cite{LS95}, corresponds to solving the semilinear equation
\[
-\Delta \mathrm{U}(x) = \mathrm{f}(|x|, \mathrm{U}(x)) \quad \text{ in } \quad  \mathbb{R}^3,
\]
which, in the framework where
$$
\mathrm{f}(|x|, \mathrm{U}(x)) = |x|^\alpha |\mathrm{U}(x)|^{p-2}\mathrm{U}(x),\quad  \alpha > 0,\,\,\, p > 2,
$$
becomes the now well-known H\'{e}non equation:
\[
-\Delta \mathrm{U}(x) = |x|^\alpha |\mathrm{U}(x)|^{p-2}\mathrm{U}(x) \quad \text{ in } \quad  \mathbb{R}^3.
\]
We must observe that in the scenario of astrophysics, the weight \( \mathfrak{h}(x) = |x|^\alpha \) represents a black hole situated at the cluster's center,
with its absorption intensity increasing as \( \alpha \) grows.

Currently, elliptic models of H\'{e}non-type serve as standard prototypes for several diffusion phenomena. From a mathematical perspective, the problem
\begin{equation}\label{EqHenon}
\left\{
\begin{array}{rclcl}
-\Delta u(x) & = & |x|^{\alpha} |u|^{p-2} u & \text{in } & \Omega, \\
u(x) & = & 0 & \text{on } & \partial \Omega,
\end{array}
\right.
\end{equation}
has a fascinating structure, and numerous results have been established. As an example, Ni in \cite{Ni82}, in the context of problem \eqref{EqHenon}, noted that the inclusion of the weight function $\mathfrak{h}(x) = |x|^{\alpha}$ affects the Poho\u{z}aev identity, leading to a novel critical exponent, namely $p_{\alpha} = \frac{2(n + \alpha)}{n - 2}$, for the existence of classical solutions. Finally, in \cite{dosSP16} (see also references therein), the authors studied the concentration profiles of various types of symmetric positive solutions.

Additionally, for a fixed $p>1$ (and $u > 0$), the model equation \eqref{EqHenon} is referred to as a Hardy-type equation if $\alpha < 0$ because of its connection to the Hardy-Sobolev inequality (cf. \cite{CotLa19} for a survey), and as a H\'{e}non-type equation if $\alpha > 0$. Finally, for $\alpha \in \mathbb{R}$, this model is known as the Hardy-H\'{e}non equation. In these contexts, the following results are available in the literature:

By using variational techniques, one can establish the existence of a nontrivial solution to \eqref{EqHenon} in $H_0^1(\Omega)$ (for $\Omega \subset \mathbb{R}^n$ an arbitrary bounded domain) provided that
    $$
    1 < p < \frac{2(n - |\alpha|)}{n-2}-1,
    $$
    through an application of the Caffarelli-Kohn-Nirenberg estimates (see \cite{CKN84}). Moreover, by employing a generalized Pohozaev-type identity, one proves the non-existence of nontrivial solutions in star-shaped domains if
    $$
    0 \geq \alpha > -n \quad \text{and} \quad p = \frac{2(n - |\alpha|)}{n-2}-1.
    $$

We emphasize that the simplest model equation described by \eqref{EqHenon} (for non-negative weak solutions) serves as our initial impetus for studying Hardy-H\'{e}non-type equations for more general quasi-linear models.

\subsection{Matukuma and Batt– Faltenbacher–Horst type equations}\label{Subsection1.2}

Consider the quasilinear elliptic equation
\begin{equation}\label{Eq6.1}
\mathrm{div}\left(\mathfrak{a}(|x|) |\nabla u|^{p-2} \nabla u\right) = \mathfrak{h}(|x|) f(u) \quad \text{in } \Omega,
\end{equation}
where $ 1 < p < \infty$, and \(\mathfrak{a}, \mathfrak{h} : \mathbb{R}^+ \to \mathbb{R}^+\) are radial functions that are sufficiently regular, \(\Omega\) is a domain of \(\mathbb{R}^n\), \(n \geq 1\), containing the origin, and $f:  \mathbb{R}_{+} \to \mathbb{R}_{+}$ is a non-decreasing function (see, Remark \ref{Remark-M-BFH}).

As discussed by Pucci-Serrin in \cite{PucSer06}, remember that prototypes of \eqref{Eq6.1}, with nontrivial functions \(\mathfrak{a}\), \(\mathfrak{h}\), arise, for example, from equations of Matukuma type and Batt–Faltenbacher–Horst type (cf. \cite{BFH86}). Specifically, the Matukuma-type equation is given by
\begin{equation}\label{Eq6.2}
\Delta_p u = \frac{f(u)}{1 + r^\sigma}, \quad r = |x|, \quad \text{for} \quad \sigma > 0,
\end{equation}
where \(p > 1\), \(\mathfrak{a}(|x|) \equiv 1\), and \(\mathfrak{h}(|x|) = \frac{1}{1 + r^\sigma}\).

A second model example is the equation
\begin{equation}\label{Eq6.3}
\Delta_p u = \frac{r^\sigma}{(1 + r^{p^{\prime}})^{\frac{\sigma}{p^{\prime}}}} \cdot \frac{f(u)}{r^{p^{\prime}}} \quad \text{for} \quad \sigma > 0,
\end{equation}
where now \(\mathfrak{a}(|x|) \equiv 1\) and \(\mathfrak{h}(|x|) = \frac{r^{\sigma - p^{\prime}}}{(1 + r^{p^{\prime}})^{\frac{\sigma}{p^{\prime}}}}\) with $r=|x|$.

It is important to highlight that the Batt–Faltenbacher–Horst equation, introduced in \cite{BFH86}, arises in astrophysics as a model of stellar structure. This equation is obtained when $p = 2$ and $n = 3$, and it reduces to the standard Matukuma equation when $\sigma = 2$.

All these equations are discussed in detail in \cite[section 4]{PG-HMS06}, as special cases of the main example
\begin{equation}\label{Eq6.4}
\mathrm{div}\left(r^k | \nabla u|^{p-2} \nabla u\right) = r^{l} \left(\frac{r^s}{1 + r^s}\right)^{\frac{\sigma}{s}} f(u), \quad k \in \mathbb{R}, \, l \in \mathbb{R}, \, s > 0 \quad \text{and} \quad \sigma > 0,
\end{equation}
where now \(\mathfrak{a}(|x|)  = r^k \) and  \( \mathfrak{h}(|x|) = r^{l} \left(\frac{r^s}{1 + r^s}\right)^{\frac{\sigma}{s}}\), with $r=|x|$.

Particularly, in \cite{PG-HMS06}, conditions on the exponents were identified, ensuring that, under appropriate behavior of the nonlinearity \(f\), radial ground states for \eqref{Eq6.2}-\eqref{Eq6.4} are unique.

Finally, in \cite{FilPucRig08} (see also \cite{FilPucRig09}), the authors provided sufficient conditions for the non-existence of non-negative, non-trivial entire weak solutions of class \(W^{1,p}_{\text{loc}}(\mathbb{R}^n) \cap C(\mathbb{R}^n)\) for \(p\)-Laplacian elliptic inequalities with possibly singular weights
$$
\mathrm{div}\left(\mathfrak{a}(|x|) |\nabla u|^{p-2} \nabla u\right) \geq  \mathfrak{h}(|x|) f(u) \quad \text{in } \mathbb{R}^n\setminus \{0\} \quad \text{with} \quad u\geq 0.
$$
This inequality arises in the study of non-Newtonian fluids, non-Newtonian filtration, and subsonic gas motion (cf. \cite{Yang} and references therein), and also in Riemannian geometry (cf. \cite{PRS05}). In conclusion, the importance of such inequalities has recently been widely recognized (see \cite{PG-HMS06} and references therein).

In conclusion, our regularity estimates immediately result in optimal regularity for Matukuma and Batt–Faltenbacher–Horst type models with general weights, which may be singular (see Remark \ref{Remark-M-BFH} and Example \ref{Example6.7}). In this direction, we believe that our estimates can also play a significant role in certain geometry contexts.

\subsection{Quasilinear models and its regularity theories}

For \( 1< p<  \infty \), consider the \( p \)-Laplace equation
\begin{equation}\label{Eq-p-Laplace}
  -\Delta_p u \defeq - \text{div} \left( |\nabla u|^{p-2} \nabla u \right) = 0 \quad \text{in} \quad \Omega,
\end{equation}
where \( \Omega \subset \mathbb{R}^n \) is a bounded and open domain. Recall that weak solutions to the quasilinear equation \eqref{Eq-p-Laplace} are widely known as $p$-harmonic functions. Note that for \( p > 2 \), the model equation \eqref{Eq-p-Laplace} is degenerate elliptic, and for \( 1 < p < 2 \), it is singular, since the modulus of ellipticity $\mathbf{a}(\nabla u) = |\nabla u|^{p-2}$ degenerates (or blows up) at those points where \( |\nabla u| = 0 \).

Regarding regularity estimates, in the late sixties, Ural’tseva in \cite{Ural68} proved that for $p > 2$, weak solutions of equation \eqref{Eq-p-Laplace} have H\"{o}lder continuous first derivatives in the interior of $\Omega$. Almost one decade later, Uhlenbeck in \cite{Uhl77} proved a far-reaching extension to a certain class of nonlinear elliptic systems. In the early eighties, Lewis in \cite{Lewis83} and DiBenedetto in \cite{DiBe} gave proofs valid for the singular scenario, i.e., $1 < p < 2$. Nevertheless, in general, weak solutions do not have any better regularity than $C_{\text{loc}}^{1, \alpha}$ for some unknown $\alpha  \in (0, 1)$.

In \cite{DiBe}, by considering suitable structural assumptions on vector fields \(\mathbf{a}\) and \(\mathbf{b}\), Evans proved \(C^{1,\alpha}\)-regularity of weak solutions of elliptic equations of the form
\[
-\operatorname{div}(\mathbf{a}(x, u, \nabla u)) + \mathbf{b}(x, u, \nabla u) = 0.
\]
In particular, such assumptions cover model equations of the $p$-Laplacian type as \eqref{Eq-p-Laplace}. Furthermore, these results apply also to inhomogeneous equations with general forcing terms. We should highlight that similar results were derived independently by Tolksdorf in \cite{Tolks}.

In the last decade, Kuusi and Mingione in \cite{KuusiMing12} derived pointwise bounds for weak solutions to potentially degenerate equations of the form
\[
-\text{div} \, \mathfrak{A}(x, \nabla u) = \mu \quad \text{in} \, B_1,
\]
where \( \mathfrak{A} : \mathbb{R}^n \times B_1 \to \mathbb{R}^d \) satisfies the standard structural conditions, and \( \mu : B_1 \to \mathbb{R} \) is a Borel measure with finite total variation. The conditions imposed on \( \mathfrak{A} \) are flexible enough to allow for the presence of coefficients. In this direction, by using Wolff-type potentials associated with \( \mu \), the authors reveal a mechanism that transfers regularity properties from the data to the solutions. As a result, they establish both Schauder and Calder\'{o}n-Zygmund estimates, in H\"{o}lder and Sobolev spaces, respectively (cf. Duzar-Mingione's work \cite{DuzMing11} and Mingione's work \cite{Ming11} for enlightening related topics).

At this point, we must compare our estimates from Theorems \ref{sharp_regularity} and \ref{Hessian_continuity} with the ones obtained by Teixeira in \cite[Theorem 3]{Tei14}. In that work, Teixeira addressed the sharp gradient regularity, on critical sets, of solutions to quasilinear equations of the form
\begin{equation}\label{TeixPDE}
\div(\mathfrak{a}(x, \nabla u)) = \mu \quad \text{in} \quad \Omega \subset \mathbb{R}^n, \,\,\,n \ge 2,
\end{equation}
where $\mu$ is a source function or a measure with finite total mass, i.e., $|\mu|(\Omega)< \infty$ and such that $\mu(x) \in L^q(\Omega)$ $(q>n)$.

In this scenario, under the assumptions \eqref{condestr}-\eqref{intdoscoef}, Teixeira establishes \( u \in C_{\text{loc}}^{1, \min\left\{\alpha_{\mathrm{H}}^{-}, \frac{q-n}{q(p-1)}\right\}} \), on the set of critical points \( \mathcal{S}_u(\Omega):= \{x \in \Omega: |\nabla u(x)|=0\} \), where \( \alpha_{\mathrm{H}} \in (0, 1] \) denotes the maximal H\"{o}lder exponent of an \( \mathfrak{a} \)-harmonic profile, i.e., a weak solution of the equation
$$
\div(\mathfrak{a}(\nabla \mathfrak{h})) = 0 \quad \text{in} \quad \Omega.
$$

In particular, when \( f \in L^q(B_1)\cap C^0(B_1) \), our estimates improve upon those provided by Teixeira in \cite[Theorem 3]{Tei14}.

\subsubsection{Quasilinear equations: recent trends in regularity theory}

From now on, let us consider the following quasilinear elliptic equation in divergence form
\begin{equation}\label{ToyModel}
  \mathrm{div} \,\mathfrak{a}(x, \nabla u) = f(x, u) \quad \text{in} \ \Omega \subset \mathbb{R}^n,
\end{equation}
for a suitable vector field $\mathfrak{a} \in C^0(\Omega \times \mathbb{R}^n; \mathbb{R}^n)$ and an appropriate source term $f: \Omega \times \mathbb{R} \to \mathbb{R}$ under growth or integrability conditions to be clarified soon (see \eqref{condestr}-\eqref{EqHomog-f}).

In the context of elliptic regularity theory, an intriguing and challenging issue that arises is determining the optimal conditions on the inhomogeneous source term $f(x, s)$, along with appropriate constraints on the leading coefficients $x \mapsto \mathfrak{a}(x, \xi)$, and possibly dependence on the spatial dimension $n$, in order to achieve the sharp modulus of continuity for the gradient of weak solutions.

By way of motivation, in the framework of the semilinear elliptic equation
\begin{equation}\label{EqSemilin}
  \Delta u = f(x, u) \quad \text{in} \quad  B_1,
\end{equation}
Shahgholian in \cite[Theorem 1.2]{Shah03} established the interior $C^{1,1}$ regularity of $u$ under the assumptions that

$$
 \textbf{Assumption} \,\,\,{\bf \mathrm{H}}:
 \left\{
 \begin{array}{l}
   x \mapsto f(x, u) \quad \text{is Lipschitz, uniformly in } u; \\
   \partial_u f \geq -\mathrm{C} \,\text{ weakly for some } \mathrm{C} \in \mathbb{R}.
 \end{array}
 \right.
 $$

Nevertheless, the precise conditions on $f(x, u)$ that guarantee the $C^{1,1}$ regularity of $u$ were provided by Indrei \textit{et al.} \cite[Theorem 1.1]{IMN17}. Precisely, they demonstrated the interior $C^{1,1}$ regularity of $u$ under the following structural conditions:
 $$
 \textbf{Assumption} \,\,\,{\bf \mathrm{A}}:
 \left\{
 \begin{array}{l}
   \mbox{Let } f = f(x, u) \,\mbox{ be Dini continuous in } u, \,\, \mbox{ uniformly in } x;\\ \\
   \mbox{Assume further } f \,\mbox{ has a } C^{1,1} \,\mbox{ Newtonian potential in } x, \,\, \mbox{ uniformly in } u.
 \end{array}
 \right.
 $$
In some configurations, such conditions, namely in the \textbf{Assumption} {\bf $\mathrm{A}$}, are sharp. Moreover, such assumptions embrace functions that fail to fulfill both conditions in Shahgholian’s theorem.

Now, concerning the quasilinear scenario, recently the authors in \cite{ATU17} examined the $p$-Poisson equation
\begin{equation}\label{EqPlane}
  -\Delta_p u = f(x) \quad \text{in} \quad  B_1
\end{equation}
for $p > 2$, where they assumed that $f \in L^\infty(B_1)$ and $B_1 \subset \mathbb{R}^n$ is the $n$-dimensional unity ball. The central result of this manuscript ensures that for any bounded weak solution $u \in W^{1,p}(B_1)$ of \eqref{EqPlane}, that $u \in C^{1, \min\left\{\alpha^{-}_{\mathrm{M}}, \frac{1}{p-1}\right\}}(B_{1/2})$. Furthermore,
\[
\|u\|_{C^{1, \min\left\{\alpha^{-}_{\mathrm{M}}, \frac{1}{p-1}\right\}}(B_{1/2})} \leq \mathrm{C}_p \cdot\left(\|u\|_{L^{\infty}(B_1)} +  \|f\|_{L^\infty(B_1)}^{\frac{1}{p-1}} \right),
\]
where $\mathrm{C}_p > 0$ is a constant depending only on $p$ and dimension.

In certain problems, obtaining sharp (potentially enhanced) control over the behavior of solutions in specific regions or at particular points becomes crucial for advancing research in a wide class of nonlinear elliptic models. By way of illustration, this is a central concern in the theory of free boundary problems, e.g., quasilinear dead core problems (see \cite{daSRosSal19}, \cite{daSSal18}, and \cite{Diaz85})
$$
\Delta_p\, u = \lambda u_{+}^q \quad \text{in} \quad B_1 \quad \text{for} \quad \lambda>0, \,0<q<p-1 \quad \text{and} \quad p \in (1, \infty),
$$
as well as one-phase Bernoulli-type problems driven by the $p$-Laplacian (see \cite{DanPetros05})
$$
\Delta_p\, u = 0 \quad \text{in} \quad B_1 \cap \{u>0\},\quad u=0, \,\,|\nabla u| = \frac{1}{(p-1)^{\frac{1}{p}}}\quad \text{on} \quad \partial\{u>0\}\cap B_1,
$$
or even in various problems with a geometric appeal, e.g., obstacle-type problems (see \cite{ALS15} and \cite{BJrDaST23}), namely finding the smallest solution to the problem:
$$
\left\{
\begin{array}{rcl}
\div\left( |\nabla u|^{p-2}\mathfrak{A}(x) \nabla u\right)  = f(x) & \text{in} & \{u> \phi\}\cap B_1,\\
\div\left( |\nabla u|^{p-2}\mathfrak{A}(x) \nabla u\right)  \le f(x) & \text{in} & B_1,\\
u(x)\ge \phi(x) & \text{in} & B_1,
\end{array}
\right.
$$
where
\begin{equation}\label{Cond_Obst_Prob}
|\mathfrak{A}(x) - \mathfrak{A}(y)|\leq |x-y|^{\kappa}, \,(\kappa \in (0, 1]) \quad \phi \in C^{1, \beta}(B_1) \,\,(\beta \in (0, 1]) \quad \text{and} \quad f \in L^q(B_1) \text{ with } q \geq \frac{np}{p-1}.
\end{equation}

Finally, as a final motivation, we must mention Teixeira's work in \cite{Tei22}, where regularity estimates at interior stationary points of solutions to $p$-degenerate elliptic equations in an inhomogeneous medium were studied
$$
\mathrm{div}\,\mathfrak{a}(x, \nabla u) = f(x, u) \lesssim \mathrm{c}_0|x|^{\alpha}|u|^{m}
$$
for some $\mathrm{c}_0>0$, $\alpha \geq 0$, and $0 \leq m < p - 1$ (for $p > 2$). In such a scenario, the data satisfy the structural conditions \eqref{condestr}-\eqref{EqDini-Cond}.

In this context, if the source term is bounded away from zero, Teixeira derives a quantitative non-degeneracy estimate, indicating that solutions cannot be smoother than $C^{p^{\prime}}$ at stationary points (see \cite[Proposition 5]{Tei22}). Furthermore, at critical points where the source vanishes, Teixeira establishes higher-order regularity estimates that are sharp regarding the vanishing rate of the source term (see \cite[Theorem 3]{Tei22}).

Subsequently, we summarize the sharp regularity estimates reported in the literature for several problems related to \eqref{ToyModel} in the table below:

\begin{table}[h]
\centering
\resizebox{\textwidth}{!}{
 \begin{tabular}{c|c|c|c}
{\it Toy model } & {\it Structural conditions }  & {\it Sharp regularity estimates} & \textit{References} \\
\hline
 $\Delta u = f(x, u)$ & $\textbf{Assumption}$ \,\,\,${\bf \mathrm{H}}$ & $C_
{\text{loc}}^{1, 1}$ &  \cite{Shah03}  \\\hline
  $\Delta u = f(x, u)$ & $\textbf{Assumption}$ \,\,\,${\bf \mathrm{A}}$ & $C_
{\text{loc}}^{1, 1}$ &  \cite{IMN17}  \\\hline
 $-\Delta_p u=f(x)$ & $p>2$ \text{and} $f \in L^{\infty}$ & $C_{\text{loc}}^{1, \min\left\{\alpha^{-}_{\mathrm{H}}, \frac{1}{p-1}\right\}}$ &  \cite{ATU18} \\
 \hline
  $ -\div(\mathfrak{a}(x, \nabla u)) = \mu$ & \eqref{condestr}-\eqref{EqDini-Cond} \text{and}\,\,$p>2$ & $C_{\text{loc}}^{1, \min\left\{\alpha_{\mathrm{H}}^{-}, \frac{q-n}{q(p-1)}\right\}}$ &  \cite{Tei14} \\
 \hline
 $\Delta_p\, u = \lambda u_{+}^q(x)  $ & $\lambda>0 , \,q \in [0, p-1) \quad \text{and} \quad p \in (1, \infty)$ & $C_
{\text{loc}}^{\frac{p}{p-1-q}}$ & \cite{daSRosSal19}, \cite{daSSal18} and \cite{Diaz85}  \\
\hline
$\div\left( |\nabla u|^{p-2}\mathfrak{A}(x) \nabla u\right)  = f(x)  \,\,\text{in} \,\, \{u> \phi\}$ & \eqref{Cond_Obst_Prob} $\quad \text{and} \quad p \in (1, \infty)$ & $C_
{\text{loc}}^{1, \min\left\{\beta,\min\left\{ \kappa,1-\frac{n}{q}\right\}\cdot \min\left\{1, \frac{1}{p-1}\right\}\right\}}$ & \cite{ALS15} and \cite{BJrDaST23} \\
\hline
$\Delta_p u = f(x, u) \lesssim \mathrm{c}_0|x|^{\alpha}|u|^{m}$ & $\alpha \geq 0$ and $0\le m<p-1$ (for $p>2$) & $C_
{\text{loc}}^{1, \min\left\{\alpha_{\mathrm{H}}, \frac{p+\alpha+1}{p-1-m}\right\}^{-}}$ & \cite{Tei22} \\
\end{tabular}}
\end{table}

In particular, we should highlight that our findings will improve some results previously presented, to some extent, via alternative strategies and techniques.


\subsection*{Quasilinear models in the planar scenario: Quasiregular mappings}

Additionally, in addition to the previously cited contributions, we must stress the sharp regularity in dimension two. More specifically, Aronsson in \cite{Arons89} and Iwaniec-Manfredi in \cite{IwaManf89} proved that a \( p \)-harmonic function belongs to \( C^{k,\alpha}_{\text{loc}} \cap W_{\text{loc}}^{2+k, q} \), where the integer \( k \geq 1 \) and the exponent \( \alpha \in (0, 1) \) is given by the equation
\begin{equation}\label{SharpExp}
k + \alpha^{\sharp}_p = \frac{1}{6} \left( 7 + \frac{1}{p - 1} + \sqrt{1 + \frac{14}{p - 1} + \frac{1}{(p - 1)^2}} \right) \quad \text{and}  \quad 1 \leq q(p) < \frac{2}{2-\alpha}.
\end{equation}

The proof substantially exploits the ideas from Bojarski and Iwaniec's fundamental survey \cite{BojIwa84}. In this context, a key ingredient is the hodograph method, which transforms a \( p \)-harmonic equation into a linear first-order elliptic system. This system is solved using a Fourier series method. Furthermore, a careful examination of the Fourier expansion formula for the system's solution leads us to the desired regularity statement. In particular, they demonstrate that the regularity exponent in \eqref{SharpExp} is optimal. We also refer interested readers to Manfredi's work \cite{Manf88} and Aronsson-Lindqvist's manuscript \cite{AronLind88}, which are closely related to Bojarski and Iwaniec's work \cite{BojIwa84}.

Recently, in the seminal paper \cite{ATU17}, the authors considered the \( p \)-Poisson equation
\begin{equation}\label{Eqp-Poisson}
  -\Delta_p u = f(x), \quad \text{in} \quad B_1 \subseteq \mathbb{R}^2
\end{equation}
for \( p > 2 \), with the additional assumption that \( f \in L^\infty(B_1) \). Thus, they addressed a longstanding conjecture: the proof of the \( C^{p^{\prime}} \)-regularity conjecture in the plane, where this regularity is optimal:
\begin{theorem}[{\cite{ATU17}}]
Let \( B_1 \subset \mathbb{R}^2 \), and let \( u \in W^{1,p}(B_1) \) be a weak solution of
\[
-\Delta_p u = f(x), \quad p > 2,
\]
with \( f \in L^\infty(B_1) \). Then, \( u \in C^{p^{\prime}}(B_{1/2}) \) and
\[
\|u\|_{C^{p^{\prime}}(B_{1/2})} \leq \mathrm{C}_p \left( \|u\|_{L^{\infty}(B_1)} + \|f\|^{\frac{1}{p-1}}_{L^\infty(B_1)}\right).
\]
\end{theorem}

Note that in the above statement, the constant \( p^{\prime} \) is the well-known H\"{o}lder conjugate of \( p \), i.e.,
$$
p + p^{\prime} = pp^{\prime} \quad \Leftrightarrow \quad p^{\prime} = \frac{1}{1-\frac{1}{p}} = 1 + \frac{1}{p-1}.
$$

It is important to note that the regularity estimates are a consequence of a new oscillation-type estimate. Specifically, the authors establish that there exists \( \delta_0 > 0 \) such that if \( \|f\|_{L^{\infty}(B_1)} \leq \delta_0 \) and \( u \in W^{1,p}(B_1) \) is a weak solution of \eqref{Eqp-Poisson} with \( \|u\|_{L^{\infty}(B_1)} \leq 1 \), then there exists a constant \( \mathrm{C}(p) > 1 \) for which
\[
\sup_{x \in B_r} |u(x) - u(0)| \lesssim \mathrm{C}(p) \left(r^{p^{\prime}} + |\nabla u(0)| r \right) \quad \text{for all} \quad r > 0.
\]

Another crucial ingredient in \cite{ATU17} is that the conjecture holds true provided that \( p \)-harmonic functions, i.e., weak solutions of \eqref{Eq-p-Laplace}, are locally of class \( C^{1, \alpha_p} \) for some \( \alpha_p > \frac{1}{p - 1} \). While this remains a challenging and unresolved problem in higher dimensions (cf. [3]), it holds true in the plane, thereby providing a complete account of this conjecture in 2-D. The crucial estimate follows from Baernstein II and Kovalev's pivotal result in \cite{BaerKova05}, which exploits the fact that the complex gradient of a \( p \)-harmonic function in the plane (i.e., \( \phi = \frac{\partial u}{\partial z} \)) turns out to be a \( \mathbf{K} \)-quasiregular gradient mapping, with \( \mathbf{K} = \mathbf{K}(p) = p-1 \). The complex gradient \( \phi = u_x - i u_y \in W^{1, 2}(B_1) \) satisfies the system
\[
\phi_{\overline{z}} = \left(\frac{1}{p}-\frac{1}{2}\right)\left(\frac{\phi_z \overline{\phi}}{\phi} + \frac{\overline{\phi}_z \phi}{\overline{\phi}}\right), \quad \text{with} \quad z = x + iy,
\]
in the hodograph plane. Therefore, by using the known H\"{o}lder exponents for quasiregular mappings (see \cite[Section 5]{BaerKova05} and \cite{LL}), it follows that a weak solution of the \( p \)-Laplacian is of class \( C_{\text{loc}}^{1, \alpha^{\ast}_p} \) with
$$
 \alpha^{\ast}_p  =  \frac{1}{2p}\left(-3-\frac{1}{p-1}+ \sqrt{33+\frac{30}{p-1} + \frac{1}{(p-1)^2}}\right).
$$
Furthermore, it holds
$$
[\phi]_{C^{0, \alpha^{\ast}_p}(B_{1/2})} \leq \mathrm{C}_p\|u\|_{L^{\infty}(B_1)}.
$$

Next, we will summarize this result in the following proposition.

\begin{proposition}[{\cite{ATU17} and \cite{BaerKova05}}]\label{PropB-K} For any \( p > 2 \), there exists \( 0 < \tau_0 < \frac{p-2}{p-1} \) such that \( p \)-harmonic functions in \( B_1 \subset \mathbb{R}^2 \) are locally of class \( C^{p^{\prime} + \tau_0} \). Furthermore, if \( u \in W^{1,p}(B_1) \cap C^0(B_1) \) is \( p \)-harmonic in the unit disk \( B_1 \subset \mathbb{R}^2 \), then there exists a constant \( \mathrm{C}_p>0 \), depending only on \( p \), such that
\[
[\nabla u]_{C^{0, \frac{1}{p-1} + \tau_0}(B_{1/2})} \leq \mathrm{C}_p \| u \|_{L^\infty(B_1)}.
\]
\end{proposition}

Now, it is worth emphasizing that, we have
$$
\alpha_p^{\sharp} > \alpha^{\ast}_p> \frac{1}{p-1} \quad \text{for any}  \quad p > 2.
$$

Finally, we would like to refer interested readers to Astala \textit{et al.}'s work \cite{ACFJK20}, which addresses improved H\"{o}lder regularity for strongly elliptic PDEs. Specifically, they deal with the interaction of two elliptic PDEs in two dimensions, namely the Leray-Lions equation
$$
\text{div}\mathcal{A}(z, \nabla u) = 0,
$$
and the nonlinear Beltrami equation
$$
\phi_{\overline{z}} = \mathcal{H}(z, \phi_z).
$$

The structural field \( \mathcal{H}(z, \zeta) \) is assumed to be measurable in \( z \in \Omega \), and ellipticity is quantified by assuming uniform Lipschitz control
\begin{equation}\label{eqUnifEllip}
|\mathcal{H}(z, \zeta) - \mathcal{H}(z, \eta)| \leq \kappa |\zeta - \eta|, \quad \zeta, \eta \in \mathbb{C},
\end{equation}
where \( 0 \leq \kappa < 1 \) is a fixed constant.

In the case of autonomous equations,
\begin{equation}\label{EqAuton}
  \phi_{\overline{z}} = \mathcal{H}(\phi_z),
\end{equation}
the solutions are shown to exhibit H\"{o}lder regularity that is superior to the classical one. Indeed, it is well-known that the solutions of \eqref{EqAuton} have $\mathrm{K}$-quasiregular directional derivatives with
$\mathrm{K} = \frac{1+\kappa}{1-\kappa}$ (see \cite{BojIwa84} and \cite{HM20}). For this reason, by invoking Morrey's classical embedding theorem (see \cite[Section 3.10]{AIM09}), the solutions belong to $C_{\text{loc}}^{1, \frac{1}{\mathrm{K}}}$.

\begin{theorem}[{\cite[Theorem 1]{ACFJK20}}]\label{ThmACFJK} Under the ellipticity assumption \eqref{eqUnifEllip}, solutions \( \phi \) to the autonomous Beltrami equation \eqref{EqAuton}
belong to \( C^{1,\alpha_{\mathrm{K}}}_{\text{loc}} \), where for \( \mathrm{K} = \frac{1+\kappa}{1-\kappa}\),
\[
\alpha_{\mathrm{K}} = \frac{1 - \kappa}{1 + \frac{\kappa}{2}\max \{ 1, 2 - 4\kappa \}} = \min \left\{ \frac{4}{3\mathrm{K} + 1}, \frac{\mathrm{K} + 1}{3\mathrm{K} - 1} \right\} > \frac{1}{\mathrm{K}}.
\]
\end{theorem}

\begin{remark}
Similar to the results presented in Corollary \ref{Corol02}, we can consider Leray-Lions type equations in the planar setting,
$$ \text{div}\, \mathcal{A}(\nabla u) = 0, $$
which, as discussed above, is closely related to the nonlinear autonomous Beltrami equation \eqref{EqAuton}. Moreover, under the assumption \( \alpha_{\mathrm{K}} \in \left(\frac{1+\alpha}{p-1}, 1\right) \), it implies that \( u \in C^{1, \frac{\alpha+1}{p-1}}(B_{1/2}) \). In addition, \( u \) satisfies the estimate
\[
\|u\|_{C^{1, \frac{\alpha+1}{p-1}}(B_{1/2})} \leq \mathrm{C}_{p, n, \alpha} \left( \|u\|_{L^{\infty}(B_1)} + \|f\|^{\frac{1}{p-1}}_{L^\infty(B_1)}\right).
\]
\end{remark}

\subsection{Preliminaries and auxiliary results}

To begin with, we will recall the definition of weak solutions:

\begin{definition}[{\bf Weak solutions}]
We say that a function $u$ is a weak solution of
$$
-\div \,\mathfrak{a}(x,\nabla u) = f(x) \quad \text{in} \quad \Omega
$$
if $u \in W^{1,p}_{loc}(\Omega)$ and 	
$$
\int_{\Omega} \mathfrak{a}(x,\nabla u) \cdot \nabla \Psi \, dx = \int_{\Omega} f\Psi \, dx
$$
for every test function $\Psi \in W^{1,p}_0(\Omega)$. \\

In a similar manner, we define weak sub-solutions and weak super-solutions. We say that a function $u$ is a weak sub-solution (resp. weak super-solution) of
$$
-\div \,\mathfrak{a}(x,\nabla u) = f(x) \quad \text{in} \quad \Omega
$$
if $u \in W^{1,p}_{loc}(\Omega)$ and
\begin{eqnarray*}
\int_{\Omega} \mathfrak{a}(x,\nabla u) \cdot \nabla \Psi \, dx \geq \int_{\Omega} f\Psi \, dx \quad \left(\text{resp.} \quad \int_{\Omega} \mathfrak{a}(x,\nabla u) \cdot \nabla \Psi \, dx \leq \int_{\Omega} f\Psi \, dx \right)
\end{eqnarray*}
for every test function $\Psi \in W^{1,p}_0(\Omega)$, where $\Psi \geq 0$.
\end{definition}

\begin{theorem}[{\bf Comparison Principle \cite[Theorem 2.4.1]{pucciserrin}}]\label{comp}
Let $u$ e $v$ be, respectively, weak solutions to
$$
\div \,\mathfrak{a}(x,\nabla v) +\mathrm{B}(x, v) \le 0 \leq   \div \,\mathfrak{a}(x,\nabla u) +\mathrm{B}(x, u) \quad  \text{in} \quad  \Omega.
$$
Suppose that $\mathfrak{a}: \Omega \times \mathbb{R}^n \to \mathbb{R}^n$ satisfies the monotonicity condition
\begin{eqnarray*}
\langle \mathfrak{a}(x,\xi) -\mathfrak{a}(x, \eta),\xi - \eta \rangle  > 0 \quad \text{when} \quad \xi \neq \eta.
\end{eqnarray*}
while $\mathrm{B} = \mathrm{B}(x, z)$ is non-increasing in $z$. If $u\leq v$ on $\partial \Omega$, then $u\leq v$ in $\Omega$.
\end{theorem}

Next, we present Campanato's Embedding Theorem, which is useful in obtaining the estimates from Theorem \ref{sharp_regularity}. Hereafter, we denote
$$
(f)_{x_0,r} := \dfrac{1}{|B_r(x_0)|} \int_{B_r(x_0)} f(x) \, dx = \intav{B_r(x_0)} f(x) \, dx.
$$

\begin{theorem}[{\bf Campanato's Theorem \cite[Theorem 1.54]{MZ97}}] \label{campanato}
Let $u \in L^2(\Omega)$ and $0 < \alpha \leq 1$. Suppose there exists a constant $\mathbf{M}_0 > 0$ such that
$$
\intav{B_r(x_0)} |u - u_{x_0,r}|^2 \, dx \leq \mathbf{M}_0 r^{2\alpha}
$$
for any ball $B_r(x_0) \subset \Omega$. Then, $u \in C^{0,\alpha}(\Omega)$, and for any ball $B_r(x_0) \subset \Omega$, the following holds:
$$
\operatorname{osc}_{B_{r/2}(x_0)} u \leq \mathrm{C} \mathbf{M}_0 r^{\alpha},
$$
where $\mathrm{C} = \mathrm{C}(n, \alpha)$.
\end{theorem}

Next, suppose $\mathfrak{a}_{ij} \in L^{\infty}(B_1)$ is a uniformly elliptic matrix in $B_1 = B_1(0)$, i.e.,
\[
\lambda |\xi|^2 \leq \mathfrak{a}_{ij}(x) \xi_i \xi_j \leq \Lambda |\xi|^2 \quad \text{for any } x \in B_1, \, \xi \in \mathbb{R}^n.
\]
Now, we assume that $u \in H^1(B_1)$ satisfies the weak formulation
\begin{equation}\label{EqWeaksolLinear}
\int_{B_1} \mathfrak{a}_{ij} D_i u D_j \varphi \, dx = \int_{B_1} f \varphi \, dx \quad \text{for any } \varphi \in H^1_0(B_1).
\end{equation}
Then, the following H\"{o}lder estimates for gradients hold true.

\begin{theorem}[{\cite[Theorem 3.13]{QhanFlin}}]\label{casep=2} Let $u \in H^1(B_1)$ be a weak solution to \eqref{EqWeaksolLinear}. Assume $\mathfrak{a}_{ij} \in C^{0, \alpha}(\overline{B_1})$, and $f \in L^q(B_1)$ for some $q > n$ and $\alpha = 1 - \frac{n}{q} \in (0, 1)$. Then, $\nabla u \in C^{0, \alpha}(B_1)$.
Moreover, there exists an $R_0 = R_0(\lambda, \|a_{ij}\|_{C^{0, \alpha}})$ such that for any $x \in B_{1/2}$
and $r \leq R_0$, there holds
\[
\displaystyle \intav{B_r(x)}
|\nabla u - (\nabla u)_{x,r}|^2 \leq \mathrm{C} r^{2\alpha} \left\{ \|f\|_{L^q(B_1)}^2 + \|u\|_{H^1(B_1)}^2 \right\}
\]
where $\mathrm{C} = \mathrm{C}(\lambda, \|\mathfrak{a}_{ij}\|_{C^{0, \alpha}})$ is a positive constant.
\end{theorem}

\begin{theorem}[{\bf Stability of solutions}]\label{estabilidade}(\cite[Theorem 2.1, Remark 2.1, Remark 2.5]{Almostconv})
Assume that
\begin{equation}\label{1}
|\mathfrak{a}(x,\xi)|  \leq \mathrm{k}_2|\xi|^{p-1}
\end{equation}
\begin{equation}\label{22}
\langle \mathfrak{a}(x,\xi)-\mathfrak{a}(x,\eta), \xi-\eta\rangle  >0
\end{equation}
\begin{equation}\label{33}
\dfrac{\langle \mathfrak{a}(x,\xi),\xi\rangle}{|\xi|+|\xi|^{p-1}} \to \infty \quad  \text{when} \quad |\xi| \to \infty
\end{equation}
for almost every $x\in \Omega$ and for all $\xi,\eta \in \mathbb{R}^n, \xi\neq \eta$. Consider the non linear elliptic equations in the weak sense
\begin{eqnarray}\label{44}
\left\{
\begin{array}{rclr}
- \	\div \, \mathfrak{a}(x,\nabla u_n)& = & f_n(x) & \text{in } \  \Omega \\
u_n(x) & = & 0 & \text{on} \  \partial  \Omega
\end{array}
\right.
\end{eqnarray}
and assume that
\begin{equation}\label{55}
u_n \rightharpoonup u  \quad \text{in} \quad W^{1,p}(\Omega), \quad u_n \rightarrow u \quad \text{in} \quad L^p_{loc}(\Omega) \quad \text{and} \quad u_n \rightarrow u \quad \text{a.e. in} \quad \Omega	
\end{equation}
and
\begin{equation}\label{6}
f_n \rightarrow f \quad \text{in} \quad (W^{1,p}_{0})^{\prime}.
\end{equation}

Then, we can passing the limit in \eqref{4}, this is, $u$ is a weak solution to
\begin{eqnarray}
\left\{
\begin{array}{rclr}
- \	\div \, \mathfrak{a}(x,\nabla u)& = & f(x) & \text{in } \ \Omega \\
u(x) & = & 0 & \text{on} \  \partial  \Omega.
\end{array}
\right.
\end{eqnarray}
\end{theorem}
\medskip

Next, we can state the main regularity estimate for inhomogeneous quasilinear elliptic models with varying coefficients.

\begin{theorem}[{\cite[Theorem 4]{DuzMing10}, \cite[Theorem 1.6]{KuusiMing12} and \cite[Theorem 5.10]{DZ23}}]\label{ModUnivEst} Let $u \in W^{1, p}(\Omega)$ be a weak solution to \eqref{EqModel}, where $\mathfrak{a}$ satisfies the structural conditions \eqref{condestr} and \eqref{EqDini-Cond}. Assume further that $f \in L^q(\Omega)$, $q>n$. Then, $u \in C^1_{\text{loc}}(\Omega)$. Moreover, for any sub-domain $\Omega^{\prime} \Subset \Omega$, there exists a universal modulus of continuity, $\tau$, depending only on $\Omega^{\prime}, \Omega, \lambda, \Lambda, \mathfrak{L}_0, \omega$ and $\|f\|_{L^q(\Omega)}$, such that
$$
\displaystyle \sup_{x, y \in \Omega^{\prime} \atop{x \neq y, |\kappa|=1}} \frac{|D^{\kappa}u(x)-D^{\kappa}u(y)|}{\tau(|x-y|)} < \infty.
$$
\end{theorem}

In order to establish the existence and uniqueness of solutions and obtain a priori estimates, universal control over the $L^{\infty}$ norm of weak solutions to an equation of the form \eqref{TeixPDE} (with $\mu = f$) is necessary. To achieve this, it is crucial to employ a version of the Aleksandrov-Bakel'man-Pucci Maximum Principle, as addressed by the authors in \cite{BJrDaST23-ABP}. In the case of bounded domains, this bound yields the following theorem.

\begin{theorem}[{\bf A.B.P. estimates  \cite[Theorem 2.4]{BJrDaST23-ABP}}]\label{abp}
Let $u \in W^{1, p}(\Omega) \cap C(\overline{\Omega})$ be a weak sub-solution of \eqref{EqModel} (resp. weak super-solution of \eqref{EqModel}). Suppose the hypotheses {\bf (H1)-(H7)} are in force.
Then, there exists a constant $\mathrm{C} > 0$, depending on $p, q, n$, and $\Lambda$, such that 	
\begin{eqnarray*}
\sup_{\Omega} u(x) \leq \sup _{\partial \Omega} u^{+}(x) + \mathrm{C} \cdot |\Omega|^{\frac{pq-n}{nq(p-1)}} \cdot \left\|f^{+}\right\|_{L^{q}\left(\Omega\right)}^{\frac{1}{p-1}}	
\end{eqnarray*}
\begin{eqnarray*}
\left(\text{resp.} \quad \inf _{\Omega} u(x) \geqslant - \inf _{\partial \Omega} u^{-}(x) - \mathrm{C} \cdot  |\Omega|^{\frac{pq-n}{nq(p-1)}} \cdot \left\|f^{-}\right\|_{L^{q}\left( \Omega\right)}^{\frac{1}{p-1}}\right).
\end{eqnarray*}

In particular, if $u$ is a weak solution of \eqref{EqModel}, then
$$
\|u\|_{L^{\infty}(\Omega)} \leq \|u\|_{L^{\infty}(\partial \Omega)} + \mathrm{C} \cdot |\Omega|^{\frac{pq-n}{nq(p-1)}} \cdot \left\|f\right\|_{L^{q}\left(\Omega\right)}^{\frac{1}{p-1}}.
$$
\end{theorem}

In conclusion, it is worth noting that, in addition to the A.B.P. estimates, Theorems \ref{comp} and \ref{estabilidade} are crucial tools for obtaining the existence and uniqueness of solutions to \eqref{pobst}.

\section{Proof of Theorem \ref{sharp_regularity}}
We begin by proving the following approximation lemma.

\begin{lemma}[\textbf{Approximation lemma}]\label{Approximation lemma}
Let $u \in W^{1,p}(B_1)$ be a weak solution to \eqref{pobst}, with $\|u\|_{L^\infty(B_1)} \leq 1$. Assume that $\mathfrak{a}$ satisfies the structural conditions \eqref{condestr}.
Given $\delta > 0$, there exists $\varepsilon > 0$, depending only on $n, \mathcal{L}_0, \omega, p, \lambda, \Lambda$, such that if
$$
|\mathfrak{a}(x,\xi) - \mathfrak{a}(0,\xi)| \leq \varepsilon |\xi|^{p-1} \quad \text{and} \quad \|f\|_{\Xi(B_1)} \leq \varepsilon,
$$
then there exists a function $\mathbf{h} \in \mathbb{H}_{\Lambda, \lambda, n}\left(B_{\frac{1}{2}}\right)$ such that
$$
\|u - \mathbf{h}\|_{L^\infty(B_1)} + \|\nabla(u - \mathbf{h})\|_{L^\infty(B_1)} \leq \delta.
$$
\end{lemma}

\begin{proof}
Assume the result does not hold. We can then find $\delta_0 > 0$ and sequences of functions $\{u_k\}_{k \in \mathbb{N}} \in L^\infty(B_1)$ and $\{f_k\}_{k \in \mathbb{N}} \in \Xi(B_1)$ such that
\begin{equation}\label{Ap_1}
\div \, \mathfrak{a}_k (x, \nabla u_k) = f_k \,\, \text{in} \,\, B_1, \quad \text{with} \quad \|u_k\|_{L^\infty(B_1)} \leq 1,\,\,\,\text{and} \,\,\, \|f_k\|_{\Xi(B_1)} = \text{o}(1),
\end{equation}
where $\mathfrak{a}_k$ satisfies \eqref{condestr} and
\begin{equation}\label{Ap_2}
|\mathfrak{a}_k(x, \xi) - \mathfrak{a}_k(0, \xi)| = \text{o}(1) |\xi|^{p-1},
\end{equation}
however
\begin{equation}\label{HC}
\|u_k - \mathbf{h}\|_{L^\infty(B_1)} + \|\nabla(u_k - \mathbf{h})\|_{L^\infty(B_1)} > \delta_0,
\end{equation}
for all functions $\mathbf{h} \in \mathbb{H}_{\Lambda, \lambda, n}\left(B_{1/2}\right)$. Now, by the $C^1$ local regularity estimates for quasilinear equations (see, Theorem \ref{ModUnivEst}), we can extract a subsequence such that
$$
u_k \rightarrow u_\infty \quad \text{and} \quad \nabla u_k \rightarrow \nabla u_\infty
$$
locally uniformly. Also, by the Arzelà–Ascoli theorem, there exists a subsequence such that $\mathfrak{a}_{k_j}(x, \cdot) \rightarrow \overline{\mathfrak{a}}(0, \cdot)$ locally uniformly. Hence, from \eqref{Ap_2}, for any $x \in B_{1/2}$ and $\xi \in B_R$, for arbitrary fixed $R > 0$, we have
$$
|\mathfrak{a}_{k_j}(x, \xi) - \overline{\mathfrak{a}}(0, \xi)| = |\mathfrak{a}_{k_j}(x, \xi) - \mathfrak{a}_{k_j}(0, \xi)| + |\mathfrak{a}_{k_j}(0, \xi) - \overline{\mathfrak{a}}(0, \xi)| = \text{o}(1).
$$
This verifies that
$$
\mathfrak{a}_{k_j}(x, \xi) \rightarrow \overline{\mathfrak{a}}(0, \xi)  = \overline{\mathfrak{a}}(\xi) \quad \text{locally uniformly in} \,\, B_{1/2} \times \mathbb{R}^n.
$$
Moreover, by passing to the limit in the PDE (see, Theorem \ref{estabilidade}), we obtain
\begin{equation}\label{u_infty}
\div(\overline{\mathfrak{a}} \, (\nabla u_{\infty})) = 0 \quad \text{in} \quad B_{1/2}, \quad \text{with} \quad \|u_{\infty}\|_{L^\infty(B_1)} \leq 1.
\end{equation}
Finally, by comparing the conclusions obtained from \eqref{u_infty} with \eqref{HC}, we reach a contradiction for $k \gg 1$ if we consider $\mathbf{h} = u_{\infty}$. The lemma is proven.
\end{proof}

\begin{lemma}\label{control_gradient}
Under the hypotheses of Lemma \ref{Approximation lemma}, there exist $\varepsilon_\star > 0$ and $0 < \rho < 1/2$ depending only on $\delta, n, \omega, p, \lambda, \Lambda$, and $\beta$ such that if
$$
|\mathfrak{a}(x, \xi) - \mathfrak{a}(0, \xi)| \leq \varepsilon_\star |\xi|^{p-1} \quad \text{and} \quad \|f\|_{\Xi(B_1)} \leq \varepsilon_\star,
$$
then the following holds
\begin{equation}
\sup_{B_\rho} |u(x) - u(0) - \nabla u(0) \cdot x| \leq \rho^{1+\beta},
\end{equation}
for $\beta$ as in \eqref{sharp_exponent}.
\end{lemma}

\begin{proof}
For a $\delta > 0$ to be chosen \textit{a posteriori}, we apply Lemma \ref{Approximation lemma} to find $\varepsilon_\star = \varepsilon_\star(\delta) > 0$ and a function $\mathbf{h}: B_{1/2} \rightarrow \mathbb{R}$ satisfying
$$
\div(\overline{\mathfrak{a}}\, (\nabla \mathbf{h})) = 0 \quad \text{in} \quad B_{{1}/{2}}
$$
in the distributional sense, and
\begin{equation}\label{Discrete_1}
    \|u - \mathbf{h}\|_{C^1(B_{1/2})} \leq \delta.
\end{equation}
From $C^{1,\alpha_{\mathrm{H}}}$ regularity theory for constant coefficient equations, there exists a universal constant $\mathrm{C} > 0$ such that
\begin{equation}\label{h-estimate}
\sup_{B_\rho}|\mathbf{h}(x) - \mathbf{h}(0) - \nabla \mathbf{h}(0) \cdot x| \leq \rho^{1+\alpha_{\mathrm{H}}},
\end{equation}
for all $\rho \in (0, 1/2)$. Now, using \eqref{Discrete_1}, \eqref{h-estimate}, and the triangle inequality, we obtain
\begin{align*}
\sup_{B_\rho}|u(x) - u(0) - \nabla u(0) \cdot x| & \leq \sup_{B_\rho}|u(x) - \mathbf{h}(x)| + \sup_{B_\rho}|\mathbf{h}(0) - u(0)| \\[0.2cm]
& \quad + \sup_{B_\rho}|\mathbf{h}(x) - \mathbf{h}(0) - \nabla \mathbf{h}(0) \cdot x| \\[0.2cm]
& \quad + \sup_{B_\rho}|\nabla(u-\mathbf{h})(0) \cdot x| \\[0.2cm]
& \leq 3\delta + \mathrm{C}\rho^{1+\alpha_{\mathrm{H}}}.
\end{align*}
Finally, since $\beta$ is fixed as in \eqref{sharp_exponent}, we choose
$$
\rho \leq \left(\frac{1}{2\mathrm{C}}\right)^{\frac{1}{\alpha_{\mathrm{H}} - \beta}} \quad \text{and} \quad \delta = \frac{1}{6}\rho^{1+\beta},
$$
and obtain that
$$
\sup_{B_\rho}|u(x) - u(0) - \nabla u(0) \cdot x| \leq \rho^{1+\beta},
$$
thus completing the proof.
\end{proof}

\begin{corollary}[\textit{1st} \textbf{step of induction}]\label{1stcase}
Under the assumptions of the previous lemma, it holds that
$$
\sup_{B_\rho}|u(x)-u(0)|\leq \rho^{1+\beta} + |\nabla u(0)|\rho.
$$
\end{corollary}

\begin{lemma}[\textit{kth} \textbf{step of induction}]\label{kthstep}
Under the assumptions of Lemma \ref{control_gradient}, one has
\begin{equation}\label{kthcase}
\sup_{B_{\rho^k}}|u(x)-u(0)|\leq \rho^{k(1+\beta)} + |\nabla u(0)|\sum_{j=1}^{k-1}\rho^{k+j\beta}.
\end{equation}
\end{lemma}

\begin{proof}
The proof will be carried out via an induction argument. The case \(k=1\) is precisely the statement of Corollary \ref{1stcase}. Now, suppose that \eqref{kthcase} holds for all values of \(l=1,\ldots,k\). Our goal is to prove it for \(l=k+1\). For this purpose, define \(u_k: B_1\rightarrow \mathbb{R}\) by
$$
u_k(x)\coloneqq \frac{u(\rho^k x)-u(0)}{\mathcal{A}_k}, \quad \text{where} \quad \mathcal{A}_k\coloneqq \rho^{k(1+\beta)} + |\nabla u(0)|\sum_{j=1}^{k-1}\rho^{k+j\beta}.
$$
Note that \(u_k\) satisfies in the weak sense
$$
\div \, \mathfrak{a}_k(x,\nabla u_k) = f_k(x),
$$
where
$$
\mathfrak{a}_k(x,\xi) = \left(\frac{\mathcal{A}_k}{\rho^k}\right)^{1-p} \mathfrak{a}\left(\rho^k x, \frac{\mathcal{A}_k \xi}{\rho^k}\right) \quad \text{and} \quad f_k(x) = \frac{\rho^{kp}}{\mathcal{A}_k^{p-1}} f(\rho^k x).
$$
Observe that \(u_k(0) = 0\) and, by the induction hypothesis, we have \(\|u_k\|_{L^\infty(B_1)} \leq 1\). Moreover,
\begin{align*}
|\mathfrak{a}_k(x,\xi)-\mathfrak{a}_k(0,\xi)| & = \left(\frac{\mathcal{A}_k}{\rho^k}\right)^{1-p} \left|\mathfrak{a}\left(\rho^k x, \frac{\mathcal{A}_k}{\rho^k} \xi\right)-\mathfrak{a}\left(0, \frac{\mathcal{A}_k}{\rho^k} \xi\right)\right| \\[0.2cm]
& \leq \omega(\rho^k |x|) |\xi|^{p-1} \ll 1,
\end{align*}
and using \eqref{EqHomog-f} with \(f(x,t) = f(x)\) and \(m=0\), we obtain
$$
\|f_k\|_{\Xi(B_1)} =
\left\{
\begin{array}{cc}
\displaystyle\|f_k\|_{L^\infty(B_1)} = \frac{\rho^{kp}}{\mathcal{A}_k^{p-1}} \sup_{B_1} |f(\rho^k x)| \leq c_n \rho^{k(1+\beta(1-p))} \|f_0\|_{L^\infty(B_1)} \ll 1, \\[0.4cm]
\displaystyle \|f_k\|_{L^q(B_1)}^q = \int_{B_1} \left|\frac{\rho^{kp}}{\mathcal{A}_k^{p-1}} f(\rho^k x)\right|^q \,dx \leq \left|\frac{\rho^{kp}}{\mathcal{A}_k^{p-1}}\right|^q \int_{B_1} |f(\rho^k x)|^q \,dx.
\end{array}
\right.
$$
Here, we have used the sharp expression \eqref{sharp_exponent}. Therefore, \(u_k\), \(\mathfrak{a}_k\), and \(f_k\) satisfy the assumptions of Approximation Lemma \ref{Approximation lemma}. Hence, we can apply Corollary \ref{1stcase} to \(u_k\) and obtain
$$
\sup_{B_\rho} |u_k(x) - u_k(0)| \leq \rho^{1+\beta} + |\nabla u_k(0)| \rho,
$$
which implies
$$
\sup_{B_\rho} \left|\frac{u(\rho^k x) - u(0)}{\mathcal{A}_k}\right| \leq \rho^{1+\beta} + \frac{\rho^{k+1} |\nabla u(0)|}{\mathcal{A}_k}.
$$
Finally, by scaling back to the unit domain, we obtain
$$
\sup_{B_{\rho^{k}}} |u(x) - u(0)| \leq \rho^{(k+1)(1+\beta)} + |\nabla u(0)| \sum_{j=1}^{k} \rho^{k+1+j\beta},
$$
thereby completing the \((k+1)\)-step of induction.
\end{proof}

\begin{lemma}\label{sharp_estimate}
Suppose that the assumptions of Lemma \ref{Approximation lemma} hold. Then, there exists a universal constant $\overline{\mathrm{M}}>1$ such that, for $\rho>0$ as stated in the conclusion of Lemma \ref{control_gradient}, we have
$$
\sup_{B_r}|u(x)-u(0)|\leq \overline{\mathrm{M}} r^{1+\beta}\left(1+|\nabla u(0)|r^{-\beta}\right),
$$
for all $r\in (0,\rho)$.
\end{lemma}

\begin{proof}
Given \( r \in (0, \rho) \), we choose \( k \in \mathbb{N} \) such that \( \rho^{k+1} \leq r \leq \rho^k \). Using Lemma \ref{kthstep}, we estimate:
\begin{align*}
\sup_{B_r} \frac{|u(x) - u(0)|}{r^{1+\beta}} & \leq \frac{1}{\rho^{1+\beta}} \sup_{B_{\rho^k}} \frac{|u(x) - u(0)|}{\rho^{k(1+\beta)}} \\[0.2cm]
& \leq \frac{1}{\rho^{1+\beta}} \left( 1 + |\nabla u(0)| \rho^{-k(1+\beta)} \sum_{j=0}^{k-1} \rho^{k+j\beta} \right) \\[0.2cm]
& \leq \frac{1}{\rho^{1+\beta}} \left( 1 + |\nabla u(0)| \rho^{-k\beta} \sum_{j=0}^{k-1} \rho^{j\beta} \right) \\[0.2cm]
& \leq \frac{1}{\rho^{1+\beta}} \left( 1 + |\nabla u(0)| \rho^{-k\beta} \frac{1}{1 - \rho^\beta} \right) \\[0.2cm]
& \leq \overline{\mathrm{M}} \left( 1 + |\nabla u(0)| r^{-\beta} \right),
\end{align*}
where \( \overline{\mathrm{M}} \coloneqq \max\left\{ \frac{1}{\rho^{1+\beta}}, \frac{1}{1 - \rho^\beta} \right\} \). This concludes the proof.
\end{proof}

We conclude this section by proving Theorem \ref{sharp_regularity}. First, by normalization, we can assume that \( \|u\|_{L^\infty(B_1)} \leq 1 \).

\begin{proof}[{\bf Proof of Theorem \ref{sharp_regularity}}]

Next, we will split the proof into two cases:

\bigskip

\noindent \textbf{Case 1}: When \( |\nabla u(0)| \leq r^{\beta} \).\\

\noindent
Using Lemma \ref{sharp_estimate}, we estimate:
\begin{eqnarray*}
\sup_{x \in B_r} |u(x) - u(0) - \nabla u(0) \cdot x| & \leq & \sup_{x \in B_r} |u(x) - u(0)| + \sup_{x \in B_r} |\nabla u(0) \cdot x| \\
& \leq & \mathrm{C} r^{1+\beta}.
\end{eqnarray*}

\noindent \textbf{Case 2}: When \( |\nabla u(0)| \geq r^{\beta} \). \\

\noindent For this case, define \( \iota \coloneqq |\nabla u(0)|^{\frac{1}{\beta}} \) and let
$$
v(x) \coloneqq \frac{u(\iota x) - u(0)}{\iota^{1+\beta}}.
$$
Clearly, we have in the weak sense
$$
\text{div}\,(\tilde{\mathfrak{a}}(x, \nabla v)) = \tilde{f}(x),
$$
where
$$
\tilde{\mathfrak{a}}(x, \xi) \coloneqq \iota^{\beta(1-p)} \mathfrak{a}\left(\iota x, \iota^{\beta} \xi\right) \quad \text{and} \quad \tilde{f}(x) = \iota^{1+\beta(1-p)} f(\iota x).
$$
Moreover,
$$
v(0) = 0 \quad \text{and} \quad |\nabla v(0)| = 1.
$$
It is clear that \( \tilde{\mathfrak{a}}(x, \xi) \) satisfies the structural conditions \eqref{condestr}. Additionally, we have \( \|\tilde{f}\|_{L^\infty(B_{1/2})} \ll 1 \). From \textbf{Case 1}, we know that \( v \) is uniformly bounded in \( L^\infty(B_{1/2}) \). Therefore, given that weak solutions of \eqref{pobst} are \( C^{1, \tau_0^{\prime}} \), it follows that:
\begin{eqnarray*}
\|v\|_{C^{1, \tau_0^{\prime}}(B_{1/2})} \leq \mathrm{C}.
\end{eqnarray*}
Now, using this and the fact that \( |\nabla v(0)| = 1 \), we can find an \( r_0 = r_0(p, \omega) \) and a constant \( \mathrm{c}_0 > 0 \) such that
$$
\mathrm{c}_0 \leq |\nabla v(x)| \leq \mathrm{c}_0^{-1} \quad \text{in} \quad B_{r_0}.
$$
Therefore, \( \hat{\mathfrak{a}}(x, \xi) \) satisfies the following conditions: for every \( x, y \in B_1 \) and \( \xi, \eta \in \mathbb{R}^n \), we have
\begin{equation}
\left\{
\begin{array}{rclcl}
|\hat{\mathfrak{a}}(x, \xi)| + |\partial_{\xi}\mathfrak{a}(x, \xi)||\xi| & \leq & \mathrm{c}_0^{2-p} \Lambda |\xi| & & \\
\lambda \mathrm{c}_0^{2-p} |\eta|^2 & \le & \langle \partial_{\xi}\hat{\mathfrak{a}}(x, \xi)\eta, \eta \rangle &  & \\
\displaystyle \sup_{x, y \in B_1 \atop{x \neq y, \,\,\,|\xi| \neq 0 }} \frac{|\hat{\mathfrak{a}}(x, \xi) - \hat{\mathfrak{a}}(y, \xi)|}{\mathrm{c}_0^{2-p} |\xi|^{p-1}} & \le & \omega(|x-y|), & &  \\
\end{array}
\right.
\end{equation}
In other words, \( \tilde{\mathfrak{a}}(x, \xi) \) fulfills \eqref{condestr} for \( p = 2 \) inside \( B_{r_0} \). By Theorem \ref{casep=2}, we have \( v \in C_{loc}^{1,\alpha_0}(B_{r_0}) \). Returning to the original domain, we obtain \( u \in C^{1,\alpha_0}(B_{r_0 \iota}) \). Particularly, we obtain the following:
\begin{eqnarray*}
\sup_{B_r} |u(x) - u(0) - \nabla u(0) \cdot x| \leq \mathrm{C} r^{1+\alpha_0} \leq \mathrm{C} r^{1+\beta},
\end{eqnarray*}
for \( r \leq r_0 \iota = r_0 |\nabla u(0)|^{\frac{1}{\beta}} \).

In conclusion, we obtain the following in both cases:
\begin{eqnarray*}
\sup_{B_r} |u(x) - u(0) - \nabla u(0) \cdot x| \leq \mathrm{C} r^{1+\beta},
\end{eqnarray*}
for every \( r \in \left(0, \frac{1}{2}\right) \) such that \( r \leq r_0 |\nabla u(0)|^{\frac{1}{\beta}} \) ({\bf Case $2$}) and when \( r \geq |\nabla u(0)|^{\frac{1}{\beta}} \) ({\bf Case $1$}).

Therefore, it remains only to obtain the optimal estimate for the range
\begin{eqnarray}\label{4}
r_0 |\nabla u(0)|^{\frac{1}{\beta}} < r < |\nabla u(0)|^{\frac{1}{\beta}}.
\end{eqnarray}

Assuming that \( r \) falls into the specified interval in \eqref{4}, we then have
\begin{eqnarray*}
\sup_{B_r} |u(x) - u(0) - \nabla u(0)| \leq \sup_{B_{\iota}} |u(x) - u(0) - \nabla u(0) \cdot x|.
\end{eqnarray*}
Hence,
$$
\sup_{B_r} |u(x) - u(0) - \nabla u(0) \cdot x| \leq \mathrm{C} \iota^{1+\beta} \leq \mathrm{C} \dfrac{\iota^{1+\alpha}}{r^{1+\beta}} r^{1+\beta} \leq \dfrac{\mathrm{C}}{r_0^{1+\beta}} r^{1+\beta}.
$$

Therefore, we obtain the desired estimate for every \( r \in \left(0, \frac{1}{2}\right) \).
\end{proof}

\begin{remark}[\textbf{Normalization and Smallness Condition}]
We conclude this section by discussing the scaling properties of equations of the form \eqref{pobst}. Let us assume that $u$ solves
$$
\div\,\mathfrak{a}(x,\nabla u) = f(x,u)\quad\text{in}\quad B_1
$$
in the weak sense. For $0<A<1$ and $B>0$ arbitrary (to be chosen \textit{a posteriori}), define $v:B_1\to \mathbb{R}$ as
$$
v(x)=\frac{u(\mathrm{A}x)}{\mathrm{B}}.
$$
Direct computations show that $v$ is a weak solution to
$$
\div\, \tilde{\mathfrak{a}}(x,\nabla v)=\tilde{f}(x,v),
$$
where the modified vector field
$\tilde{\mathfrak{a}}$ is given by
$$
\tilde{\mathfrak{a}}(x,\xi)=\left(\frac{\mathrm{B}}{\mathrm{A}}\right)^{1-p}\mathfrak{a}\left(\mathrm{A}x,\frac{\mathrm{B}}{\mathrm{A}}\xi\right),
$$
and the scaled function $\tilde{f}$ is given by
$$
\tilde{f}(x,t)=\frac{\mathrm{B}^p}{\mathrm{B}^{p-1}}f\left(\mathrm{A}x,\frac{\mathrm{B}}{\mathrm{A}}t\right).
$$
Note that
$\tilde{\mathfrak{a}}$ satisfies the same structural assumptions as $\mathfrak{a}$.

Now, by choosing
$$
\mathrm{B}\coloneqq \max\{\|u\|_{L^\infty(B_1)},1\},
$$
we can assume, without loss of generality,
that the solutions are normalized, i.e., $\|u\|_{L^\infty(B_1)}\leq 1$. Now, using the third condition of \eqref{condestr}, we obtain that
\begin{eqnarray*}
\sup_{x\in B_1 \atop |\xi|\neq 0} \frac{|\tilde{\mathfrak{a}}(x,\xi)-\tilde{\mathfrak{a}}(0,\xi)|}{|\xi|^{p-1}} &=& \left(\frac{\mathrm{B}}{\mathrm{A}}\right)^{1-p}\sup_{x\in B_1 \atop |\xi|\neq 0}\frac{|\mathfrak{a}\left(\mathrm{A}x,\frac{\mathrm{B}}{\mathrm{A}}\xi\right)-\mathfrak{a}\left(0,\frac{\mathrm{B}}{\mathrm{A}}\xi\right)|}{|\xi|^{p-1}}\\
&=& \sup_{x\in B_1 \atop |\xi|\neq 0} \frac{\left|\mathfrak{a}\left(\mathrm{A}x,\frac{\mathrm{B}}{\mathrm{A}}\xi\right)-\mathfrak{a}\left(0,\frac{\mathrm{B}}{\mathrm{A}}\xi\right)\right|}{\left|\frac{\mathrm{B}}{\mathrm{A}}\xi\right|^{p-1}}\\
&\leq & \sup_{x\in B_1 \atop |\xi|\neq 0} \omega(|\mathrm{A}x|).
\end{eqnarray*}
Moreover, by \eqref{EqHomog-f}, we have that
\begin{eqnarray*}
|\tilde{f}(x,t)|&=&\left|\frac{\mathrm{A}^p}{\mathrm{B}^{p-1}}f\left(\mathrm{A}x,\frac{\mathrm{B}}{\mathrm{A}}t\right)\right|\\
&\le & \frac{\mathrm{A}^p}{\mathrm{B}^{p-1}}c_n \mathrm{A}^\alpha \left(\frac{\mathrm{B}}{\mathrm{A}}\right)^m \|f_0\|_{\infty}\\
&=& \frac{\mathrm{A}^{p+\alpha-m}}{\mathrm{B}^{p-1-m}}c_n \|f_0\|_{\infty}.
\end{eqnarray*}
Therefore, given $\varepsilon_0>0$, by choosing
$$
\mathrm{A}\leq \max\left\{\frac{1}{2}, \omega^{-1}(\varepsilon_0),\left(\frac{\mathrm{B}^{p-1-m}\varepsilon_0}{(c_n+1)(\|f_0\|_{\infty}+1)}\right)^{\frac{1}{p+\alpha-m}}\right\},
$$
we ensure that the coefficients of $\tilde{\mathfrak{a}}$, and $\tilde{f}$ satisfy the assumptions of the Approximation Lemma \ref{Approximation lemma}, thereby showing that such hypotheses are not restrictive.
\end{remark}
\section{Proof of Theorem \ref{Hessian_continuity}}

In this section, we will prove the higher regularity at local extremum points for weak solutions of \eqref{pobst}.

\begin{proof}[{\bf Proof of Theorem \ref{Hessian_continuity}}]
For simplicity, and without loss of generality, we assume $x_0=0$. By combining discrete iterative techniques with continuous reasoning (see, for instance, \cite{CKS}), it is well established that proving estimate \eqref{Higher Reg} is equivalent to verifying the existence of a constant $\mathrm{C}>0$ such that, for all $j\in \mathbb{N}$, the following holds
\begin{equation}\label{higher1}
\mathfrak{s}_{j+1}\leq \max\left\{\mathrm{C}.2^{-\hat{\beta}(j+1)}, 2^{-\hat{\beta}}\mathfrak{s}_j\right\},
\end{equation}
where
\begin{equation}\label{beta_hat}
\mathfrak{s}_j\coloneqq \sup_{B_{2^{-j}}}u\quad\text{and}\quad \hat{\beta}\coloneqq 1+\frac{\alpha+1+m}{p-1-m}.
\end{equation}
Let us suppose, for the sake of contradiction, that \eqref{higher1} fails to hold, \textit{i.e.}, that for each $k\in \mathbb{N}$, there exists $j_k\in \mathbb{N}$ such that
\begin{equation}\label{higher2}
\mathfrak{s}_{j_k+1}> \max\left\{k2^{-\hat{\beta}(j_k+1)}, 2^{-\hat{\beta}}\mathfrak{s}_{j_k}\right\}.
\end{equation}
Now, for each $k\in \mathbb{N}$, define the rescaled function $v_k: B_1\rightarrow \mathbb{R}$ by
$$
v_k(x)\coloneqq \frac{u(2^{-j_k}x)}{\mathfrak{s}_{j_k+1}}.
$$
Note that
\begin{eqnarray}
0\leq v_k(x)\leq 2^{\hat{\beta}},\label{vk1}\\
v_k(0)=0,\label{vk2}\\
\sup_{\overline{B}_{1/2}}v_k(x)=1.\label{vk3}
\end{eqnarray}
Now, for $(x,\xi)\in B_1\times \mathbb{R}^n$, we define the vector field $\mathfrak{a}^k: B_1 \times \mathbb{R}^n \to \mathbb{R}^n$ given by
$$
\mathfrak{a}^k(x,\xi)\coloneqq \left(\frac{\mathfrak{s}_{j_k+1}}{2^{-j_k}}\right)^{1-p}\mathfrak{a}\left(2^{-j_k}x,\frac{\mathfrak{s}_{j_k+1}}{2^{-j_k}}\xi\right).
$$
From this, we obtain in the weak sense that
\begin{equation}\label{v_k-solution}
\div\, \mathfrak{a}^k(x,\nabla v_k) = f_k(x,v_k),\quad\text{where}\quad
f_k(x,t)\coloneqq \frac{(2^{-j_k})^p}{\mathfrak{s}_{j_k+1}^{p-1}}f(2^{-j_k}x,\mathfrak{s}_{j_k+1}t).
\end{equation}
Now, we claim that
$$
f_k(x,v_k)\rightarrow 0,\quad\text{as}\quad k\rightarrow \infty.
$$
Indeed, using the homogeneity of $f$, \eqref{beta_hat}, and \eqref{higher2}, we obtain
\begin{eqnarray*}
|f_k(x,v_k)|&=&\frac{(2^{-j_k})^p}{\mathfrak{s}_{j_k+1}^{p-1}}|f(2^{-j_k}x, \mathfrak{s}_{j_k+1}v_k)|\\
& \leq & c_n \frac{2^{-j_k(p+\alpha)}}{\mathfrak{s}_{j_k+1}^{p-1-m}}\|f_0\|_{L^\infty(B_1)}\\
&<&\frac{\mathrm{C}_1}{k^{p-1-m}}\rightarrow 0,
\end{eqnarray*}
as $k\rightarrow \infty$. Also, using the same idea as in \cite[Lemma 2.2]{CLRT14}, we obtain that
$$
\mathfrak{a}^{k}(x,\xi)\rightarrow \overline{\mathfrak{a}}(0,\xi)
$$
locally uniformly, where $\overline{\mathfrak{a}}$ is a vector field satisfying the same structural assumptions \eqref{condestr}.

Taking into account the uniform bound of $v_k$ and \eqref{v_k-solution}, we deduce that there exist two positive constants $\tau$ and $\mathrm{C}$, independent of $k$, such that
$$
v_k\in C^{1,\tau}(B_{3/4}) \quad \text{and} \quad \|v_k\|_{C^{1,\tau}(B_{3/4})}\leq \mathrm{C}.
$$
Hence, it follows from the Arzel\`{a}-Ascoli theorem that there exists a subsequence, still denoted by $v_k$, and a function $v\in C^{1,\tau'}(B_{3/4})$ such that
$$
v_k\rightarrow v \quad\text{in}\quad C^{1,\tau'}(B_{3/4}),
$$
for any $\tau'\in (0,\tau)$. Moreover, by \eqref{vk1}-\eqref{vk3}, $v$ satisfies
\begin{equation*}
\left\{
\begin{array}{rlccc}
 \div \, \overline{\mathfrak{a}}(0,\nabla v) &= 0  & \text{in}\,\, B_{3/4}, &  v\geq 0 &  \text{in}\,\, B_{3/4},\\[0.2cm]
\displaystyle\sup_{x\in B_{1/2}}v(x)&= 1, & v(0)=0.
\end{array}
\right.
\end{equation*}
By the strong maximum principle (see \textit{e.g.} V\'{a}squez's work \cite{Vasquez}), we have necessarily $v\equiv 0$ in $B_{3/4}$, which contradicts $\displaystyle\sup_{B_{1/2}} v(x)=1$.

Finally, given $r\in (0,1)$, let $j\in \mathbb{N}$ such that $2^{-(j+1)}\leq r\leq 2^{-j}$. Then, using \eqref{higher1}, we assure the existence of a universal constant $\mathrm{C}>0$ such that
$$
\sup_{B_r} u\leq \sup_{B_{2^{-j}}} u \leq \frac{\mathrm{C}}{2^{\hat{\beta}}}r^{\hat{\beta}}.
$$
Thus, the theorem is proven.
\end{proof}

As previously mentioned, with the help of Theorem \ref{sharp_regularity}, we can prove the growth control on the gradient stated in Corollary \ref{Ric}, thereby achieving a finer gradient control for solutions of
$$
\textrm{div} \, \mathfrak{a}(x,\nabla u)=f(|x|,u) \quad \textrm{in} \quad B_1
$$
at local extremum points.

\begin{proof}[{\bf Proof of Corollary \ref{Ric}}]
Assume, without loss of generality, that $0=x_0 \in \partial\{u>0\}\cap B_{1/2}$. As previously noted, it suffices to prove the following estimate
\begin{equation}\label{Grad1}
    \mathfrak{G}_{j+1}\leq \max\left\{\mathrm{C}.2^{-\tilde{\beta}(j+1)}, 2^{-\tilde{\beta}}\mathfrak{G}_j\right\},
\end{equation}
where
\begin{equation}\label{Grad2}
\mathfrak{G}_j\coloneqq \sup_{B_{2^{-j}}}|\nabla u|\quad\text{and}\quad \tilde{\beta}\coloneqq \frac{\alpha+1+m}{p-1-m}.
\end{equation}
Suppose, for the sake of contradiction, that \eqref{Grad1} does not hold. Then, for each $k\in \mathbb{N}$, there exists $j_k\in \mathbb{N}$ such that
\begin{equation}\label{Div2}
\mathfrak{G}_{j_k+1}> \max\left\{k2^{-\tilde{\beta}(j_k+1)}, 2^{-\tilde{\beta}}\mathfrak{G}_{j_k}\right\}.
\end{equation}
Now, for each $k\in \mathbb{N}$, define the auxiliary function as
$$
\tilde{v}_k(x)=\frac{u(2^{j_k}x)}{\mathfrak{G}_{j_k+1}},\quad\text{for}\,\, x\in B_1.
$$
Hence, using \eqref{Higher Reg} and \eqref{Div2} we obtain
$$
0\leq \tilde{v}_k(x)\leq \frac{2^{\tilde{\beta}}}{\mathfrak{G}_{j_k+1}}\leq \frac{2^{\tilde{\beta}}}{k},\quad\text{for}\quad x\in B_1.
$$
Additionally,
$$
\mathfrak{G}_{\frac{1}{2}}=1\quad\text{and}\quad \mathfrak{G}_{1}\leq 2^{\tilde{\beta}}.
$$
Now, for $(x,\xi)\in B_1\times \R^n$, we define the vector field $\tilde{\mathfrak{a}}^k: B_1 \times \mathbb{R}^n \to \mathbb{R}^n$ given by
$$
\tilde{\mathfrak{a}}^k(x,\xi)\coloneqq \left(\frac{\mathfrak{G}_{j_k+1}}{2^{-j_k}}\right)^{1-p}\mathfrak{a}\left(2^{-j_k}x,\frac{\mathfrak{G}_{j_k+1}}{2^{-j_k}}\xi\right).
$$
Thus, we obtain in the weak sense that
\begin{equation}
\div\, \tilde{\mathfrak{a}}^k(x,\nabla \tilde{v}_k) = \tilde{f}_k(x,\tilde{v}_k),\quad\text{where}\quad
\tilde{f}_k(x,t)\coloneqq \frac{(2^{-j_k})^p}{\mathfrak{G}_{j_k+1}^{p-1}}f(2^{-j_k}x,\mathfrak{G}_{j_k+1}t).
\end{equation}
Now, using the homogeneity of $f$ and \eqref{Div2}, we obtain
\begin{eqnarray*}
|\tilde{f}_k(x,\tilde{v}_k)| &=& \frac{(2^{-j_k})^p}{\mathfrak{G}_{j_k+1}^{p-1}}|f(2^{-j_k}x, \mathfrak{G}_{j_k+1}\tilde{v}_k)|\\
&\leq & c_n \frac{2^{-j_k(p+\alpha)}}{\mathfrak{G}_{j_k+1}^{p-1-m}}\|f_0\|_{L^\infty(B_1)}\\
&<& \frac{2^{(p+\alpha)}}{k^{p-1-m}}.
\end{eqnarray*}
Finally, by invoking the uniform gradient estimates from \cite{Choe, DiBe, Tolks}, we obtain
\begin{eqnarray*}
    1&=&\mathfrak{G}_\frac{1}{2}\\
& \leq & \mathrm{C}(n, p)\left[ \|\tilde{v}_k\|_{L^\infty(B_1)} + \|f(x,\tilde{v}_k)\|_{L^\infty(B_1)}^{\frac{1}{p-1}}\right]\\
&\leq & \mathrm{C}\left(\frac{2^{\tilde{\beta}}}{k}+\frac{2^{(p+\alpha)}}{k^{p-1-m}}\right)\rightarrow 0,
\end{eqnarray*}
as $k\rightarrow\infty$, which yields a contradiction, thereby proving the lemma.
\end{proof}


\section{Non-degeneracy of solutions: Proof of Theorem \ref{NãoDeg} }

In this section, we will establish the non-degeneracy of solutions in some specific scenarios (cf. \cite{BJrDaST23}, \cite{CLRT14} and \cite{daSSal18}).

\begin{proof}[{\bf Proof of Theorem \ref{NãoDeg}}]

Let $x_0 \in B_1$ be a critical point of $u$. Assume, without loss of generality, that $x_0 = 0$, and consider
$$
v(x) \coloneqq \epsilon_0 |x|^{1 + \frac{1 + \alpha}{p - 1}}
$$
as a comparison function. We claim that for suitable values of $\epsilon_0>0$, there exists a universal $r^{\ast} > 0$ such that
$$
\div \, \mathfrak{a}(x, \nabla v) - \mathrm{c}_0 |x|^{\alpha} \leq 0,
$$
for all $r \in (0, r^{\ast})$, $x_0 \in \mathcal{S}_u(\Omega$, and for all $x \in B_r(x_0) \subset B_1$.

To prove this claim, we will compute $\nabla v$ and estimate $| \div \, \mathfrak{a}(x, \nabla v) |$.
\begin{eqnarray}\label{gradv}
\nabla v(x) = \epsilon_0 \left(\frac{p + \alpha}{p - 1}\right) |x|^{\frac{1 + \alpha}{p - 1} - 1} x.
\end{eqnarray}
Now, observe that
\begin{eqnarray*}
\div \, \mathfrak{a}(x, \nabla v) = \sum_{i=1}^n \frac{\partial a_i}{\partial x_i}(x, \nabla v) + \sum_{i, j=1}^n \frac{\partial a_i(x, \nabla v)}{\partial (\nabla v)_j} \cdot \frac{\partial (\nabla v)_j}{\partial x_i}(x).
\end{eqnarray*}

Next, we have
\begin{eqnarray*}
|\div \, \mathfrak{a}(x, \nabla v)| &\leq& \sum_{i=1}^n \left|\frac{\partial a_i}{\partial x_i}(x, \nabla v)\right| + \sum_{i, j=1}^n \left|\frac{\partial a_i(x, \nabla v)}{\partial (\nabla v)_j}\right| \cdot \left|\frac{\partial (\nabla v)_j(x)}{\partial x_i}\right|.
\end{eqnarray*}
Now, by using \eqref{7.1} and the structural conditions \eqref{condestr}, we have that
\begin{eqnarray*}
|\div \ \mathfrak{a}(x,\nabla v)|  \leq  \Lambda_0 \left|\nabla v\right|^{p-1}+ \Lambda|\nabla v|^{p-2} \sum_{i, j=1}^n \left| \frac{\partial^2 v(x)}{\partial x_i \partial x_j}\right|:= \mathcal{S}_1+\mathcal{S}_2.
\end{eqnarray*}
Next, we will analyze $\mathcal{S}_1$ and $\mathcal{S}_2$ separately.
\begin{eqnarray}\nonumber
\mathcal{S}_1  =  \Lambda_0\left[\epsilon_0\left(\frac{p+\alpha}{p-1}\right) |x|^{\frac{1+\alpha}{p-1}}\right]^{p-1} = \Lambda_0\epsilon_0^{p-1}\left(\frac{p+\alpha}{p-1}\right)^{p-1} |x|^{1+\alpha}
\end{eqnarray}
Now, note that
\begin{eqnarray*}
\frac{\partial v}{\partial x_j} = \epsilon_0\left(\frac{p+\alpha}{p-1}\right)|x|^{\frac{1+\alpha}{p-1}-1}x_j.
\end{eqnarray*}
Then, we have
\begin{eqnarray*}
\frac{\partial^2 v}{\partial x_i \partial x_j} = \epsilon_0\left(\frac{p+\alpha}{p-1}\right)\left(\frac{1+\alpha}{p-1}-1\right)|x|^{\frac{1+\alpha}{p-1}-3}x_ix_j + \epsilon_0\left(\frac{p+\alpha}{p-1}\right)|x|^{\frac{1+\alpha}{p-1}-1}\delta_{ij},
\end{eqnarray*}
which implies the following
\begin{eqnarray}\label{hessv}
\sum_{i,j=1}^n\left|\frac{\partial^2 v}{\partial x_i \partial x_j}\right| &\leq &
\epsilon_0\left(\frac{p+\alpha}{p-1}\right)\left(\frac{1+\alpha}{p-1}-1\right)|x|^{\frac{1+\alpha}{p-1}-3}\sum_{i,j=1}^n|x_ix_j|\nonumber\\
& & + \epsilon_0\left(\frac{p+\alpha}{p-1}\right)|x|^{\frac{1+\alpha}{p-1}-1}\sum_{i,j=1}^n\delta_{ij}\nonumber\\
&=& n\epsilon_0\left(\frac{p+\alpha}{p-1}\right)\left(\frac{1+\alpha}{p-1}-1\right)|x|^{\frac{1+\alpha}{p-1}-1}\nonumber\\
& & + n \epsilon_0\left(\frac{p+\alpha}{p-1}\right)|x|^{\frac{1+\alpha}{p-1}-1}\nonumber\\
&=& n\epsilon_0\left(\frac{p+\alpha}{p-1}\right)\left(\frac{1+\alpha}{p-1}\right)|x|^{\frac{1+\alpha}{p-1}-1}.
\end{eqnarray}
Next, we will estimate \( \mathcal{S}_{2} \). Specifically, from \eqref{gradv} and \eqref{hessv}, we have
\begin{eqnarray*}
\mathcal{S}_2 & \leq & \Lambda \left[ \epsilon_0 \left( \frac{p+\alpha}{p-1} \right) |x|^{\frac{1+\alpha}{p-1}} \right]^{p-2} \left[n \epsilon_0 \left( \frac{p+\alpha}{p-1} \right) \left( \frac{1+\alpha}{p-1} \right) |x|^{\frac{1+\alpha}{p-1}-1}\right] \\[0.2cm]
& = & n \Lambda \epsilon_0^{p-1} \left( \frac{p+\alpha}{p-1} \right)^{p-1} \left( \frac{1+\alpha}{p-1} \right) |x|^{\alpha}.
\end{eqnarray*}

Thus, we obtain
\begin{eqnarray*}
|\div \ \mathfrak{a}(x,\nabla v)| & \leq &  \Lambda_0 \epsilon_0^{p-1} \left( \frac{p+\alpha}{p-1} \right)^{p-1} |x|^{1+\alpha} \\[0.2cm]
& & + n \Lambda \epsilon_0^{p-1} \left( \frac{p+\alpha}{p-1} \right)^{p-1} \left( \frac{1+\alpha}{p-1} \right) |x|^{\alpha} \\[0.2cm]
& \leq & \epsilon_0^{p-1} \left( \frac{p+\alpha}{p-1} \right)^{p-1} \left[ \Lambda_0 + n \Lambda \left( \frac{1+\alpha}{p-1} \right) \right] |x|^{\alpha}.
\end{eqnarray*}

Thus, if we choose
$$
\epsilon_0 = \left( \frac{c_0}{ \Lambda_0 + n \Lambda \left( \frac{1+\alpha}{p-1} \right)} \right)^{\frac{1}{p-1}} \left( \frac{p-1}{p+\alpha} \right),
$$
and using \eqref{7.2}, we have shown that for universally small \( r \in (0, r^{\ast}) \), it holds that
\begin{eqnarray*}
\div \ \mathfrak{a}(x,\nabla v) - c_0 |x|^\alpha & \leq & |\div \ \mathfrak{a}(x,\nabla v)| - c_0 |x|^\alpha \\[0.2cm]
& \leq & 0 \\[0.2cm]
& \leq & \div \ \mathfrak{a}(x,\nabla u) - f(x) \\[0.2cm]
& \leq & \div \ \mathfrak{a}(x,\nabla u) - c_0 |x|^\alpha.
\end{eqnarray*}

Now, for any ball \( B_r \subset B_1 \), there must exist a point \( y_r \in \partial B_r \) such that \( u(y_r) \geq v(y_r) \). Indeed, if this were not the case, then by the Comparison Principle (see Theorem \ref{comp}), we would have \( u < v \) in the entire ball \( B_r \). However, we know that \( v(0) = 0 < u(0) \). In conclusion, we can estimate
$$
\sup_{\partial B_r} u(x) \geq u(y_r) \geq v(y_r) = \left( \frac{c_0}{ \Lambda_0 + n \Lambda \left( \frac{1+\alpha}{p-1} \right)} \right)^{\frac{1}{p-1}} \left( \frac{p-1}{p+\alpha} \right) r^{1+\frac{1+\alpha}{p-1}}.
$$
Thus, the theorem is proven.
\end{proof}


\section{Proofs of Liouville-type Theorems}
As a consequence of Theorem \ref{Hessian_continuity}, we will prove Theorem \ref{Liouville I}.

\begin{proof}[{\bf Proof of Theorem \ref{Liouville I}}]
For each \( k \in \mathbb{N} \), we define
$$
u_k(x) = \frac{u(kx)}{k^{\hat{\beta}}}, \quad \text{for } x \in B_1,
$$
where \( \hat{\beta} \) is as in \eqref{beta_hat}. Note that \( u_k(0) = 0 \) and \( u_k \) is a weak solution of
$$
\div\, \mathfrak{a}_k(x, \nabla u_k) = f_k(x, u_k),
$$
where
$$
\mathfrak{a}_k(x, \xi) \coloneqq
k^{(\beta - 1)(1 - p)} \mathfrak{a}(kx, k^{\beta - 1} \xi) \quad \text{and} \quad f_k(x, t) = k^{p - \beta(p - 1)} f(kx, k^\beta t).
$$
Note that \( \mathfrak{a}_k \) satisfies \eqref{condestr}, and \( f_k \) satisfies \eqref{EqHomog-f}. From Theorem \ref{Hessian_continuity}, we deduce that if \( x_k \in \overline{B}_r \) is such that
$$
u_k(x_k) = \sup_{\overline{B}_r} u_k,
$$
where \( r > 0 \) is small, then we have
\begin{equation}\label{Liouville_eq2}
\|u_k\|_{L^\infty(B_r)} \rightarrow 0, \quad \text{as } k \rightarrow \infty.
\end{equation}
Indeed, if \( |kx_k| \) remains bounded as \( k \rightarrow \infty \), then applying Theorem \ref{Hessian_continuity} to \( u_k \), we obtain
\begin{equation}\label{Liouville_eq3}
u_k(x_k) \leq \mathrm{C}_k |x_k|^{\tilde{\beta}},
\end{equation}
where \( \mathrm{C}_k > 0 \) and \( \mathrm{C}_k \rightarrow 0 \) as \( k \rightarrow \infty \). This implies that \( u(kx_k) \) remains bounded as \( k \rightarrow \infty \), and therefore \( u_k(x_k) \rightarrow 0 \) as \( k \rightarrow \infty \), so \eqref{Liouville_eq2} holds. The result also holds in the case when \( |kx_k| \rightarrow \infty \) as \( k \rightarrow \infty \), since then, from \eqref{Liouville_eq1}, we obtain
$$
u_k(x_k) \leq |kx_k|^{-\tilde{\beta}} k^{-\tilde{\beta}} \rightarrow 0, \quad \text{as } k \rightarrow \infty.
$$
Now, if there exists \( z_0 \in \mathbb{R}^n \) such that \( u(z_0) > 0 \), by choosing \( k \in \mathbb{N} \) large enough so that \( z_0 \in B_{kr} \) and using \eqref{Liouville_eq2} and \eqref{Liouville_eq3}, we estimate
$$
\frac{u(z_0)}{|z_0|^{\tilde{\beta}}} \leq \sup_{B_{rk}} \frac{u(x)}{|x|^{\tilde{\beta}}} = \sup_{B_{r}} \frac{u_k(x)}{|x|^{\tilde{\beta}}} \leq \frac{u(z_0)}{4 |z_0|^{\tilde{\beta}}},
$$
which finally leads us to a contradiction, completing the proof of Theorem \ref{Liouville I}.
\end{proof}

Finally, we are in a position to deliver the proof of Theorem \ref{Liouville II}.

\begin{proof}[{\bf Proof of Theorem \ref{Liouville II}}]
Fix \( 0 < r_0 < R \). We consider \( v: \overline{B_R} \rightarrow \mathbb{R} \), the solution to the boundary value problem
$$
\left\{\begin{aligned}
\div (|\nabla v|^{p-2} \nabla v) & = \left( |x| - r_0 \right)^\alpha_+ v_+^m & & \text{ in } B_R \\
v & = \sup_{\partial B_R} u & & \text{ on } \partial B_R
\end{aligned}\right.
$$
By the comparison principle (Lemma \ref{comp}), we have \( u \leq v \) in \( B_R \). It follows from hypothesis \eqref{H-Liuville II} that, by taking \( R \gg 1 \) sufficiently large,
\begin{equation}\label{H1-Liouville II}
\sup_{\partial B_R} \frac{u(x)}{|x|^{1 + \frac{1 + \alpha + m}{p - 1 - m}}} \leq \theta \cdot \tau(n, m, p, \alpha)
\end{equation}
for some \( \theta < 1 \). For \( R \gg 1 \), the solution \( v = v_R \) is given by
\begin{equation}\label{H2-Liouville II}
v(x) = \tau(n, m, p, \alpha) \left( |x| - R + \left[ \frac{\partial B_R u}{\tau(n, m, p, \alpha)} \right]^{1 + \frac{1 + \alpha + m}{p - 1 - m}} \right)_{+}^{1 + \frac{1 + \alpha + m}{p - 1 - m}}.
\end{equation}
Finally, by combining \eqref{H1-Liouville II} and \eqref{H2-Liouville II}, we obtain
$$
u(x) \leq \tau(n, m, p, \alpha) \left( |x| - \left( 1 - \theta^{1 + \frac{1 + \alpha + m}{p - 1 - m}} \right) R \right)_{+}^{1 + \frac{1 + \alpha + m}{p - 1 - m}}.
$$
Letting \( R \rightarrow \infty \), we conclude the proof of the theorem.
\end{proof}

\section{Final comments }

In this final section, we will present some examples and additional comments on quasilinear diffusion models and their relationship to the problems addressed in this manuscript.

First, we will relate to the measure of critical points in quasilinear problems:
\begin{equation}\label{Eqp-Laplacian}
-\text{div} \mathfrak{a}(\nabla u) = f(x) \quad \text{in} \quad \Omega,
\end{equation}
where \( f \in L^q_{\text{loc}}(\Omega) \), with
\begin{equation}\label{Eq-Conditions}
q = \left\{
\begin{array}{rcl}
2 & \text{if } & p \geq \frac{2n}{n+2}, \\
\left(p^{\ast}\right)^{\prime} & \text{if } & 1 < p < \frac{2n}{n+2}.
\end{array}
\right.
\end{equation}
This was first addressed by Lou in \cite[Theorem 1.1]{Lou08} in the Euclidean case and under more general assumptions in the context of anisotropic quasilinear equations, as shown in the manuscript by Antonini \textit{et al.} \cite[Proposition 1.6]{ACF23}.

\begin{theorem}\label{ThmSingular-set}
 Let \( u \in W^{1,p}(\Omega) \) be a weak solution of \eqref{Eqp-Laplacian} and assume that the condition \eqref{Eq-Conditions} is satisfied. Then,
\[
f(x) = 0 \quad \text{a.e. } x \in \{x \in \Omega; |\nabla u| = 0\}.
\]
\end{theorem}

An immediate consequence of Theorem \ref{ThmSingular-set} is the following result.

\begin{corollary}
Under the assumptions of Theorem \ref{ThmSingular-set}, if \( f(x) \neq 0 \) for almost all \( x \in \Omega \), then the Lebesgue measure of the singular set \( |\nabla u| = 0 \) is zero. In particular, for any \( \mathrm{c} \in \mathbb{R} \), the level set \( u = \mathrm{c} \) has zero measure.
\end{corollary}

\begin{remark}
The previous results ensure that in our quasilinear models with forcing terms of a certain order of nullity
$$
\div \,\mathfrak{a}(x,\nabla u)  =  f(|x|, u(x)) \quad \text{in} \quad \Omega,
$$
we should have the following inclusion
$$
\mathcal{Z}_{f} = \{x \in \Omega; f(|x|, u) = 0\} \subset \mathcal{S}_u(\Omega).
$$
For this reason, this inclusion justifies why we focus on points at the boundary that are free of non-negative solutions, in order to obtain improved regularity estimates (see, Theorem \ref{Hessian_continuity} and related results).
\end{remark}

\begin{example}
Let \( \Omega = \{ x \in \mathbb{R}^n \mid 0 < \mathrm{R}_1 < |x| < \mathrm{R}_2 \} \) be an annular domain, and let \( u \in W^{1, p} \) be a weak solution to

$$
\left\{
\begin{array}{rcrcl}
  -\text{div} \left( |\nabla u|^{p-2} \nabla u \right) & = & 1 & \text{in} & \Omega \\
  u & = & 0 & \text{on} & \partial \Omega
\end{array}
\right.
$$

Then, the set of critical points is given by
\[
\mathcal{S}_u(\Omega) = \{ x \in \Omega: |\nabla u| = 0 \} = \{ x \in \Omega \mid |x| = r \}
\]
for some \( r \in (\mathrm{R}_1, \mathrm{R}_2) \). In this scenario, we have
$$
\mathscr{H}_{\text{dim}}(\mathcal{S}_u(\Omega)) = n - 1.
$$
\end{example}

\begin{example}
Let \( 1 < p < \infty \) and \( \alpha > \max\{-1, p-2\} \). Now, consider \( u \in C^{1, \beta}(\overline{B_1}) \) as a weak solution to
$$
\left\{
\begin{array}{rclcl}
  -\Delta_p u & = & |x|^{\alpha} & \text{in} & B_1 \\
  u & > & 0 & \text{in} & B_1 \\
  u & = & 0 & \text{on} & \partial B_1,
\end{array}
\right.
$$
Then, according to \cite{DP98} and \cite{DS04}, it follows that \( u \) is a radially decreasing function, with
\[
\mathcal{S}_u(B_1) := \{ x \in B_1;  |\nabla u| = 0 \} \equiv \{ 0 \}.
\]
Therefore, by arguing as in \cite{SerTang00} and applying L'Hopital's rule, we obtain
\[
|\nabla u(x)| \lesssim |x|^{\frac{1+\alpha}{p-1}} \quad \text{and} \quad \|D^2 u(x)\| \lesssim |x|^{\frac{\alpha - p + 2}{p-1}}.
\]
In particular, \( \partial_{x_i} u \) belongs to \( C^{\frac{1+\alpha}{p-1}} \) and \( \partial_{x_i x_j} u \) belongs to \( C^{\frac{\alpha - p + 2}{p-1}} \) at the origin.
\end{example}

The next result ensures that there exists a smooth regularization of a distance function,

\begin{theorem}[{\bf Calder\'{o}n-Zygmund \cite{Cald-Zyg61} and \cite[Ch. VI.2, Theorem 2]{Stein70}}]

If \( \mathbf{F} \subset \mathbb{R}^n \) is closed, then there exists a function \( f \) such that

$$
\left\{
\begin{array}{l}
  c_1(n) \, \mathrm{dist}(x, \mathbf{F}) \leq f(x) \leq c_2(n) \, \mathrm{dist}(x, \mathbf{F}), \quad \text{for all } x \in \mathbb{R}^n, \\
  f \in C^\infty(\mathbb{R}^n \setminus \mathbf{F}),
\end{array}
\right.
$$
Additionally,
\[
\left| \frac{\partial^\kappa f}{\partial x^\kappa}(x) \right| \leq c_3(n, \alpha) \, \mathrm{dist}(x, \mathbf{F})^{1 - |\kappa|} \quad \text{for all multi-indices } \kappa = (\kappa_1, \ldots, \kappa_n) \in \mathbb{N}^n.
\]
\end{theorem}

The previous theorem allow us to consider more general weights (a sort of regular distance functions) to the problem \eqref{pobst}.

\begin{example}\label{Example6.7}
Another relevant model for \eqref{pobst} is the Matukuma and Batt–Faltenbacher–Horst's problem with general weights, and noise terms	
$$
\displaystyle	\mathrm{div}(\mathrm{dist}(x, \mathrm{F})^{k_i}|\nabla u|^{p-2} \nabla u) =  \sum_{i=1}^{l_0}  \mathrm{dist}^{\alpha_i}(x, \mathrm{F_i}) u_{+}^{m_i}(x) + g_i(|x|)\quad  \text{in}  \quad  B_1
$$
where $\mathrm{F}, \mathrm{F_i} \subset B_1$ are disjoint closed sets,  $m_i \in [0, p-1)$, $\alpha_i > k_i-1-m_i$, and
$$
\displaystyle \limsup_{|x|\to 0} \frac{g_i(|x|)}{\mathrm{dist}(x, \mathrm{F}_i)^{\kappa_i}}  = \mathfrak{L}_i \in [0, \infty) \quad \text{for some} \quad \kappa_i \ge 0.
$$

Finally, in such a context, weak solutions belong to $ C_{\text{loc}}^{\displaystyle \min_{1\le i \le l_0}\left\{\frac{p+\alpha_i-k_i}{p-1-m_i}, \frac{p+\kappa_i-k_i}{p-1}\right\}}$  along the sets $\mathrm{F} \cap \mathrm{F}_i \cap \mathcal{S}_u(B_1)$.
\end{example}

\subsection*{Acknowledgments}

J.V. da Silva has received partial support from CNPq-Brazil under Grant No. 307131/2022-0 and FAEPEX-UNICAMP 2441/23 Editais Especiais - PIND - Projetos Individuais (03/2023). G.C> Ricarte has received partial support from CNPq-Brazil under Grant No. 304239/2021-6. D.S. dos Prazeres has received partial support from CNPq-Brazil under Grant No. 305680/2022-6. G.S. S\'{a} expresses gratitude for the PDJ-CNPq-Brazil (Postdoctoral Scholarship - 174130/2023-6).

\end{document}